%% file: fourier.tex
\documentclass[11pt]{article}
\usepackage{graphicx,amssymb}
\usepackage{mathtools}
\usepackage{amsmath, amsthm, amssymb, amsfonts, latexsym}
\usepackage{float}
\usepackage{color}
\usepackage[margin=2cm]{geometry}

\usepackage{etoolbox}
\usepackage{multido}

\makeatletter
\def\@@bfil{\leaders \vrule \@height \ht\z@ \@depth \z@ \hfill}
\def\@bLfil{\@@bfil}
\def\@bRfil{\@@bfil}
\def\resetbraceratio{\gdef\@bLfil{\@@bfil}\gdef\@bRfil{\@@bfil}}
\def\setbraceratio#1#2{
  \let\@bLfil\relax
  \multido{\iA=1+1}{#1}{\gappto\@bLfil{\@@bfil}}
  \let\@bRfil\relax
  \multido{\iA=1+1}{#2}{\gappto\@bRfil{\@@bfil}}
}
\def\upbracefill{$\m@th\setbox\z@\hbox{$\braceld$}\bracelu\@bLfil\bracerd\braceld\@bRfil\braceru$}
\def\downbracefill{$\m@th\setbox\z@\hbox{$\braceld$}\braceld\@bLfil\braceru\bracelu\@bRfil\bracerd$}
\makeatother

\newtheorem{theorem}{Theorem}
\newtheorem{remark}{Remark}

\DeclareMathOperator{\sinc}{sinc}

\newcommand{\dd}{\,{\rm d}}

\numberwithin{equation}{section}

\input{macpap2}

\title{A partial Fourier transform method for a class of\\ hypoelliptic Kolmogorov equations}

\author{%
C.\ Reisinger\footnote{Mathematical Institute, University of Oxford,
Andrew Wiles Building, Woodstock Rd., Oxford OX2 6GG,
\{reisinge, suli, whitley\}{@}maths.ox.ac.uk},
E.\ S{\"u}li\footnotemark[\value{footnote}], and A.\ Whitley\footnotemark[\value{footnote}]
}

\begin{document}
\maketitle

\begin{abstract}
We consider hypoelliptic Kolmogorov equations in $n+1$ spatial dimensions, with $n\geq 1$, where the differential operator in the first $n$ spatial variables featuring in the equation is second-order elliptic, and with respect to the $(n+1)$st spatial variable the equation contains a pure transport term only and is therefore first-order hyperbolic. If the two differential operators, in the first $n$ and in the $(n+1)$st co-ordinate directions, do not commute, we benefit from hypoelliptic regularization in time, and the solution for $t>0$ is smooth even for a Dirac initial datum prescribed at $t=0$. We study specifically the case where the coefficients depend only on the first $n$ variables. In that case, a Fourier transform in the last variable and standard central finite difference approximation in the other variables can be applied for the numerical solution. We prove second-order convergence in the spatial mesh size for the model hypoelliptic equation $\frac{\partial u}{\partial t} + x \frac{\partial u}{\partial y} = \frac{\partial^2 u}{\partial x^2}$ subject to the initial condition $u(x,y,0) = \delta (x) \otimes \delta (y)$, with $(x,y) \in \mathbb{R} \times\mathbb{R}$ and $t>0$, proposed by Kolmogorov, and for an extension with $n=2$. We also demonstrate exponential convergence of an approximation of the inverse Fourier transform based on the trapezium rule. Lastly, we apply the method to a PDE arising in mathematical finance, which models the distribution of the hedging error under a mis-specified derivative pricing model.
\end{abstract}



\section{Introduction}

This paper is concerned with the numerical solution of initial-value problems for a class of partial differential equations, subject to Dirac initial datum, which have the form
\begin{alignat}{2}
\label{generalpde}
\frac{\partial u}{\partial t}  + c(x,t)\frac{\partial u}{\partial y} &= \mathcal{L} u, &&\qquad (x,y,t) \in \mathbb{R}^n\times\mathbb{R}\times (0,T], \\
u(x,y,0)&=\delta (x-x_0)\otimes \delta (y-y_0), &&\qquad (x,y) \in \mathbb{R}^n\times\mathbb{R},
\label{generaldelta}
\end{alignat}
where $n \geq 1$, $T>0$, $x_0 \in \mathbb{R}^n$, $y_0 \in \mathbb{R}$, the symbol $\otimes$ signifies the (associative) binary operation of tensor product of distributions, the drift coefficient $c=c(x,t)$ is independent of the variable $y \in \mathbb{R}$, $\nabla_x c(x,t) \neq 0$ for all $(x,t) \in \mathbb{R}^n \times (0,T]$, the elliptic differential operator $\mathcal{L}$ does not include any $y$-derivatives and its coefficients only depend on $x=(x_1,\dots,x_n) \in \mathbb{R}^n$ and $t \in [0,T]$; in other words, $\mathcal{L}$ is assumed to be of the form
%
\begin{eqnarray}
\label{Lop}
\label{Lop-gen}
\mathcal{L}  = \sum_{i,j=1}^n a_{ij}(x,t) \frac{\partial^2}{\partial x_i \partial x_j} + \sum_{i=1}^n b_{i}(x,t) \frac{\partial}{\partial x_i}
+ d(x,t),
\end{eqnarray}
%
where $a_{ij}$, $i,j=1,\dots,n$, $b_i$, $i=1,\dots n$, and $d$ are continuous functions of $(x,t) \in \mathbb{R}^n \times [0,T]$, and there exists a constant $c_0>0$ such that
\[ \sum_{i,j=1}^n a_{ij}(x,t) \xi_i \xi_j \geq c_0
|\xi|^2 \qquad \forall\, \xi=(\xi_1, \dots, \xi_n) \in \mathbb{R}^n\quad \forall\, (x,t)
\in \mathbb{R}^n \times [0,T].
\]
Hypoelliptic problems of this kind arise naturally from statistical physics, stochastic analysis, and from mathematical finance in particular, as Kolmogorov equations that describe the evolution of the probability density function of stochastic processes, and the initial datum in such problems is frequently a point source, which is modelled by a Dirac measure concentrated at a point. Such initial conditions are clearly also relevant in the construction of Green's functions.

For the definition of hypoellipticity and sufficient and necessary conditions for $C^\infty$ regularity see \cite{hormander1967hypoelliptic}.
Contemporary applications and extensions to nonlinear problems are found in
\cite{villani2006hypocoercive,dric2009hypocoercivity}.

Given the interest in this class of equations, methods have recently been put forward for their numerical approximation,
with a focus on preserving the long-time behaviour of solutions to the original equation.
In \cite{foster2014structure}, a self-similar change of variables was performed and convergence of the numerical solution to the steady state under these new variables was established; furthermore, an operator splitting scheme based on decomposing the hypoelliptic operator
into coercive and convective terms was proposed.
In contrast, in \cite{PORRETTA} asymptotic properties of standard central finite difference schemes are analyzed and decay rates of difference quotients were proved in the case of $L^2$ initial data. The analysis in those papers does not cover the case of Dirac initial data, which are important in a number of applications, and indeed the numerical experiments at the end of this section show that in the case of Dirac initial datum convergence of a finite difference approximation to the problem based on central differences is not guaranteed in the discrete maximum norm.

The semidiscrete Fourier scheme that we propose here for the numerical approximation of problem \eqref{generalpde} involves the application of a Fourier transform to (\ref{generalpde}) in the $y$-direction ($y$-FT) to reduce the dimension of the problem by transforming it into a one-parameter family of parabolic problems, and it then applies finite difference discretization in the $x$-direction to this parametrized family of parabolic problems, followed by the application of an inverse Fourier transform. The observed exponential convergence of the numerical approximation to the inverse Fourier transform reduces the computational complexity of the proposed scheme to that of a finite difference approximation of a problem that has no dependence on $y$, as long as the solution is required for a single value of $y$ only. The application of the $y$-FT avoids the use of lower-order stable upwind or semi-Lagrangian discretizations 
of the $y$-derivative. Moreover, it transforms the Dirac initial datum into a constant function in the $y$-direction, which is easier
to handle numerically. To the best of our knowledge, the numerical scheme proposed in this paper is the first provably convergent numerical method for hypoelliptic problems of this type, with Dirac initial datum.

We shall study in-depth the stylized problem
\begin{alignat}{2}
\label{toypde}
\frac{\partial u}{\partial t}+x\frac{\partial u}{\partial y}&=\frac{\partial^2 u}{\partial x^2}, &&\qquad (x,y,t) \in \mathbb{R}\times\mathbb{R}\times (0,T], \\
u(x,y,0) &= \delta (x)\otimes \delta (y), &&\qquad (x,y) \in \mathbb{R} \times\mathbb{R}.
\label{toydelta}
\end{alignat}
This is the partial differential equation originally considered by Kolmogorov in his 1934 paper \cite{KOLM34}. The significance of this simple model problem stems from the fact that it incorporates two important features: hypoellipticity and Dirac initial datum. The equation (\ref{toypde}) results
from (\ref{Lop-gen}) for $d\equiv 0$ by shifting the point $(x_0,y_0)$ at which the initial Dirac datum is concentrated to the origin, linearization of the coefficient $c$ with respect to the $y$ variable around $y=0$, translation in the $y$-direction with $c(0,0) t$, and finally freezing the coefficients $c$, $b$ and $a$ with respect to $t$. Since the final term in \eqref{Lop-gen} does not affect the ideas presented herein, for the sake of simplicity of the exposition we shall confine ourselves to the case of $d\equiv 0$.
Subject to sign change, our model equation \eqref{toypde} coincides with the one studied in \cite{PORRETTA}.

Because of the special structure of the equation \eqref{toypde}, we shall apply Fourier methods in the construction of its numerical approximation and also in the convergence analysis of the proposed numerical algorithms. We shall explore the behaviour of two numerical schemes for the solution of our model problem: a semidiscrete Fourier method and a fully-discrete Fourier method, which will be described below.

Whereas the proposed numerical techniques apply to the more general model problem \eqref{generalpde}, \eqref{generaldelta},
and, in fact, with more general probability measures as initial data than the Dirac measure considered herein,
for the sake of simplicity and clarity of the exposition the mathematical analysis of the proposed numerical methods is restricted to the
simplified model \eqref{toypde}, \eqref{toydelta}, which is a special case of the Cauchy problem \eqref{generalpde}, \eqref{generaldelta} above.
For this toy model, we shall prove the convergence of the semidiscrete Fourier method and derive expressions for the leading order terms for the global discretization error.
In particular, we shall analyze the behaviour of the error between the analytical solution and its numerical approximation and will establish the rate of convergence of
the scheme as the spatial discretization parameter $\Delta x \rightarrow 0$. It should be noted that this analysis only relates to the analytical solution of the
semidiscrete scheme (\ref{eq:SEMIDSCHEME}) and its subsequent exact $y$-FT inversion. We shall discuss these additional approximations in Section \ref{sec:fourier-sol}.

We approach the task of error analysis
by applying the inverse Fourier transform to the error between $W$ and $w$, where $W$ is the solution of the equation resulting from Fourier
transforming \eqref{toypde} with respect to $y$, discretizing with respect to $x$, and applying a discrete Fourier transform with respect to $x$, whereas $w$ is the
solution of the equation resulting from Fourier transforming \eqref{toypde} with respect to both $x$ and $y$. We then perform a wavenumber analysis of the resulting
expressions to establish convergence. This method is based on similar ideas to those in
\cite{GC} and \cite{CRAW}, where time-stepping schemes for the one-dimensional heat
equation with Dirac initial datum were analyzed. While in the cited papers the Fourier transform was used purely as a mathematical tool in the analysis of the
discretization error in the original space-time co-ordinates,  here we use a partial Fourier transform (i.e., we transform in the $y$-variable only) in the construction
of the actual numerical method and use a double-Fourier transform (i.e., the Fourier transform in both the $x$ and the $y$ variable) to quantify the error of this
approximation. An interesting feature of the present analysis is the intricate interplay between the $x$- and $y$-Fourier modes, due to the hypoellipticity of the
equation.

In order to motivate the numerical method proposed in the next section, we illustrate the smoothing and convergence properties of the central difference scheme with implicit Euler time stepping,
\begin{eqnarray*}
\frac{U_{j,k}^{n+1} - U_{j,k}^{n}}{\Delta t}
+
x_j \frac{U_{j,k+1}^{n+1} - U_{j,k-1}^{n+1}}{2 \Delta y}
= \frac{U_{j+1,k}^{n+1} - 2 U_{j,k}^{n+1} + U_{j-1,k}^{n+1}}{\Delta x^2}, \qquad \mbox{with $(j,k) \in \mathbb{Z}^2$ and $n \geq 0$.}
\end{eqnarray*}
This is the scheme studied in \cite{PORRETTA}; it is shown there in particular that, for $\ell^2$ initial data, the $\ell^2$ norms of the first-order difference quotients in the $x$- and $y$-directions decay as $t^{-1}$ and $t^{-3}$, respectively, as $t \rightarrow \infty$. Closer to the situation in the present paper,
let us consider, instead, a Dirac delta concentrated at the origin in the $(x,y)$-plane. We approximate the Dirac initial datum by
\begin{eqnarray*}
U_{j,k}^0 =
\left\{
\begin{array}{cl}
\frac{1}{\Delta x\, \Delta y} & \mbox{for $(j,k)=(0,0)$,} \\
0 & \mbox{for $(j,k) \neq (0,0)$},
\end{array}
\right.
\qquad (j,k) \in \mathbb{Z}^2,
\end{eqnarray*}
which can be viewed as the mollification of the Dirac measure concentrated at the origin through convolution, in the sense of distributions, with the scaled characteristic function $\frac{1}{\Delta x\, \Delta y} \chi_{[-\Delta x/2, \Delta x/2]\times[-\Delta y/2, \Delta y/2]}$ with unit $L^1$ norm (cf. \cite{JS}). We shall consider the problem on a sufficiently large square domain $(-L,L) \times (-L,L)$ in the $(x,y)$-plane, with zero Dirichlet boundary condition along
$(\pm L, y)$ for all $y \in [-L,L]$ and along $(x,-L)$ for all $x \in [-L,L]$; in our numerical experiment below we took $L=10$, which ensures that the Dirichlet boundary condition has negligible influence on the values of the numerical solution at the final time of interest, $T=1$ in our case, close to the centre of the square where the initial Dirac delta is concentrated.

\begin{figure}[H]
\caption{Central difference scheme with Dirac initial datum; the computational domain is $x\in [-10,10]$, $y\in [-10,10]$ (smaller plot range); the numerical
solution at $T=1$ with $n_x= n_y=n_t=400$ grid spacings in the $x$-, $y$-, and $t$-directions, respectively;
left the numerical solution, right its discrete Fourier transform, with $s$ and $p$ signifying wavenumbers corresponding to the $x$ and $y$ co-ordinate
directions, respectively. 
}
\begin{center}
\includegraphics[trim={0.5cm 1 1 1},clip, width=3.4in,height=2.2in]{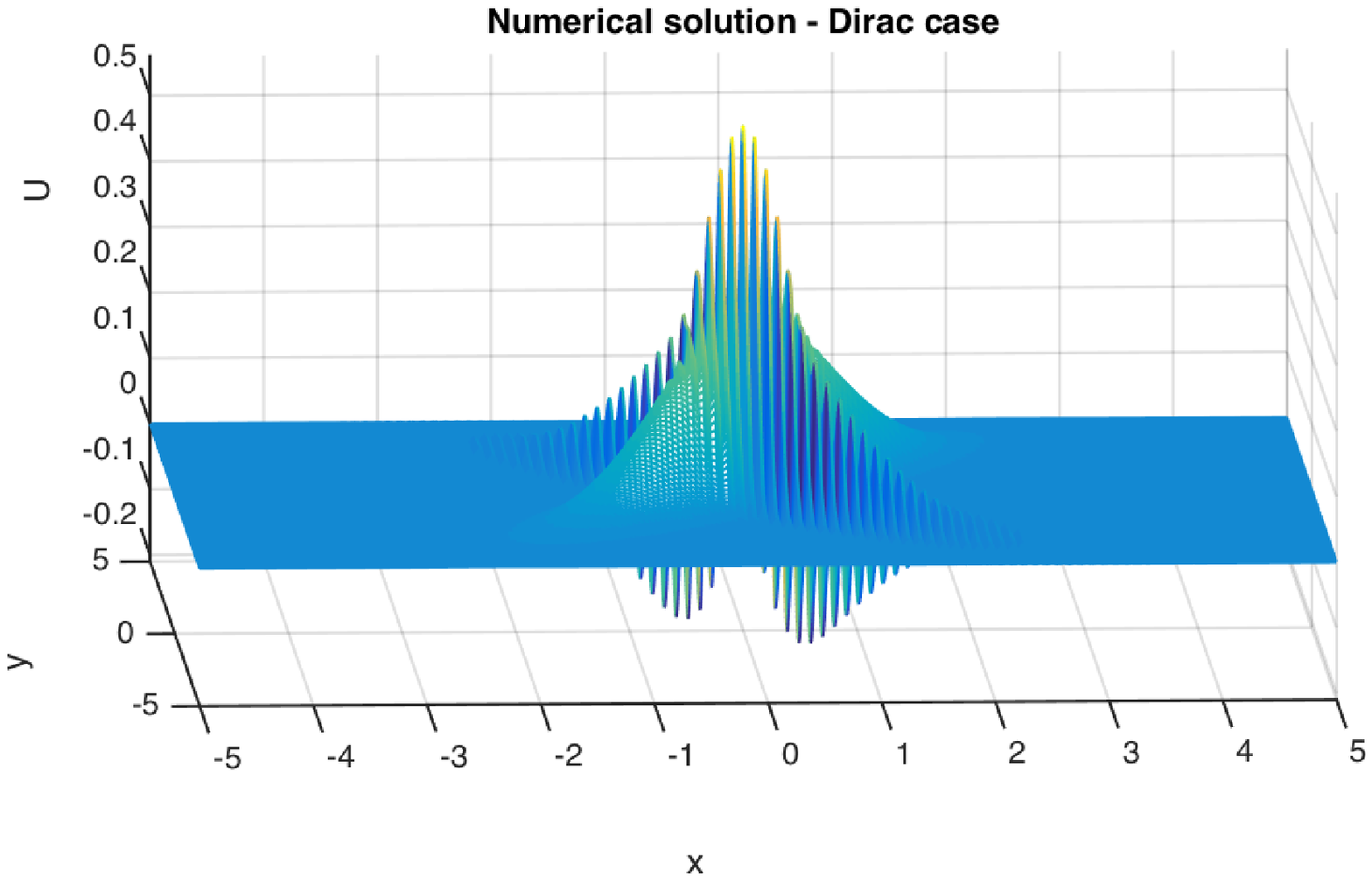} \hfill
\includegraphics[trim={0.5cm 1 1 1},clip, width=3.4in,height=2.2in]{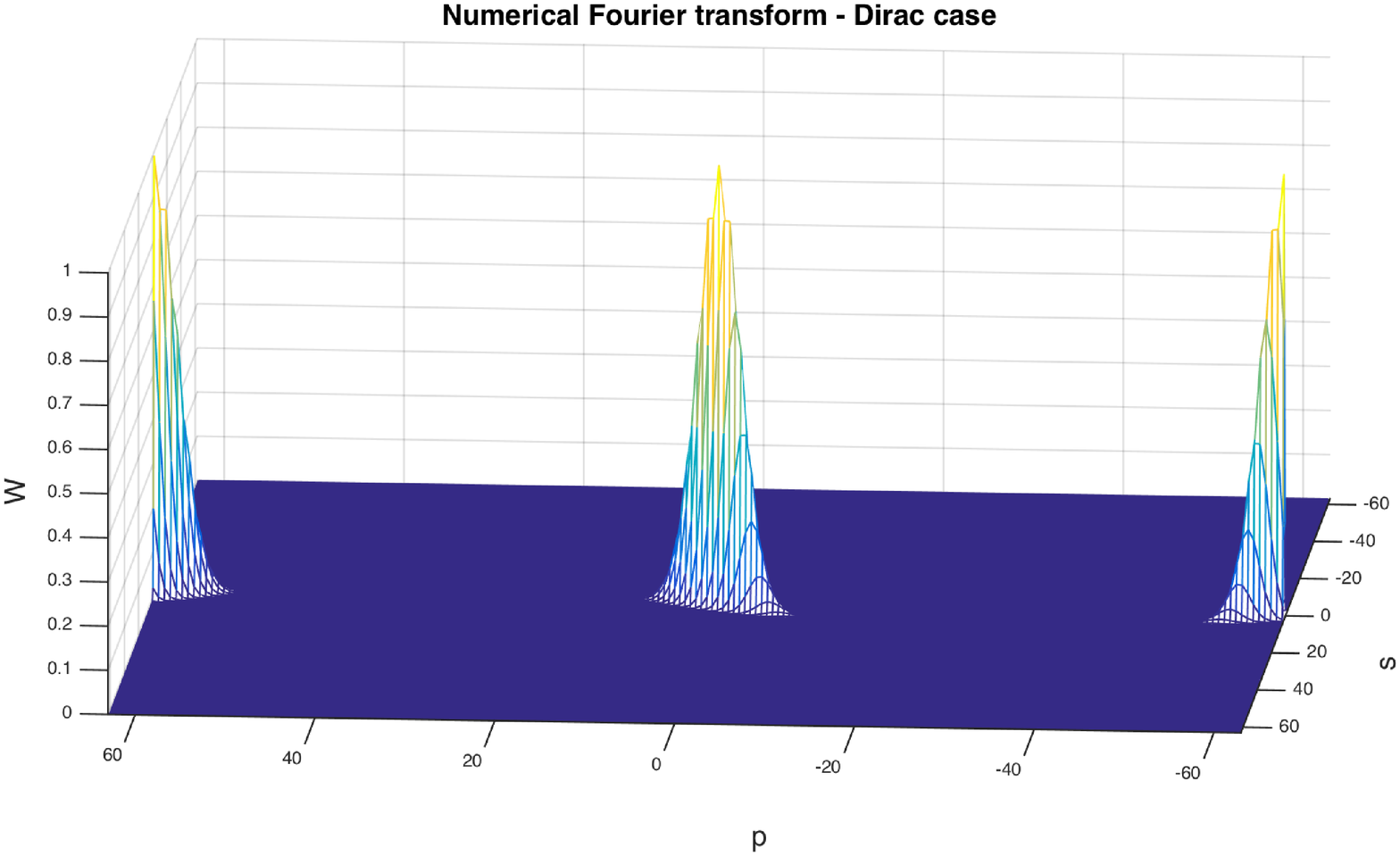}
\end{center}
\label{fig:nonsmoothcentral}
\end{figure}

Fig.~\ref{fig:nonsmoothcentral}, left, shows the numerical solution at $T=1$, which exhibits large oscillations in the $x$-direction, but is smooth in the
$y$-direction. Indeed, good approximation to the analytical solution is visible between the oscillations. The discrete Fourier transform of the numerical solution,
with wavenumbers $s$ and $p$ (described in more detail later), is depicted in the right panel in Fig.~\ref{fig:nonsmoothcentral}. It shows low wavenumber
components in both $s$ and $p$ near the origin in Fourier space, which approximate well the Fourier transform of the analytical solution, and low $s$-/high $p$-wavenumber
components concentrated at $(s,p) = (0,\pm \pi/\Delta y) = (0,\pm 20\pi) \approx (0,\pm 63)$, which trigger the spurious oscillations in the numerical approximation of the analytical solution.

\begin{figure}[H]
\caption{Central difference scheme with smooth initial datum: $x\in [-10,10]$, $y\in [-10,10]$, $t \in [1,2]$, $n_x = n_y = n_t = 400$; the numerical solution at $T=2$ (left) and
its discrete Fourier transform (right). 
}
\begin{center}
\includegraphics[trim={0.5cm 1 1 1},clip, width=3.4in,height=2.2in]{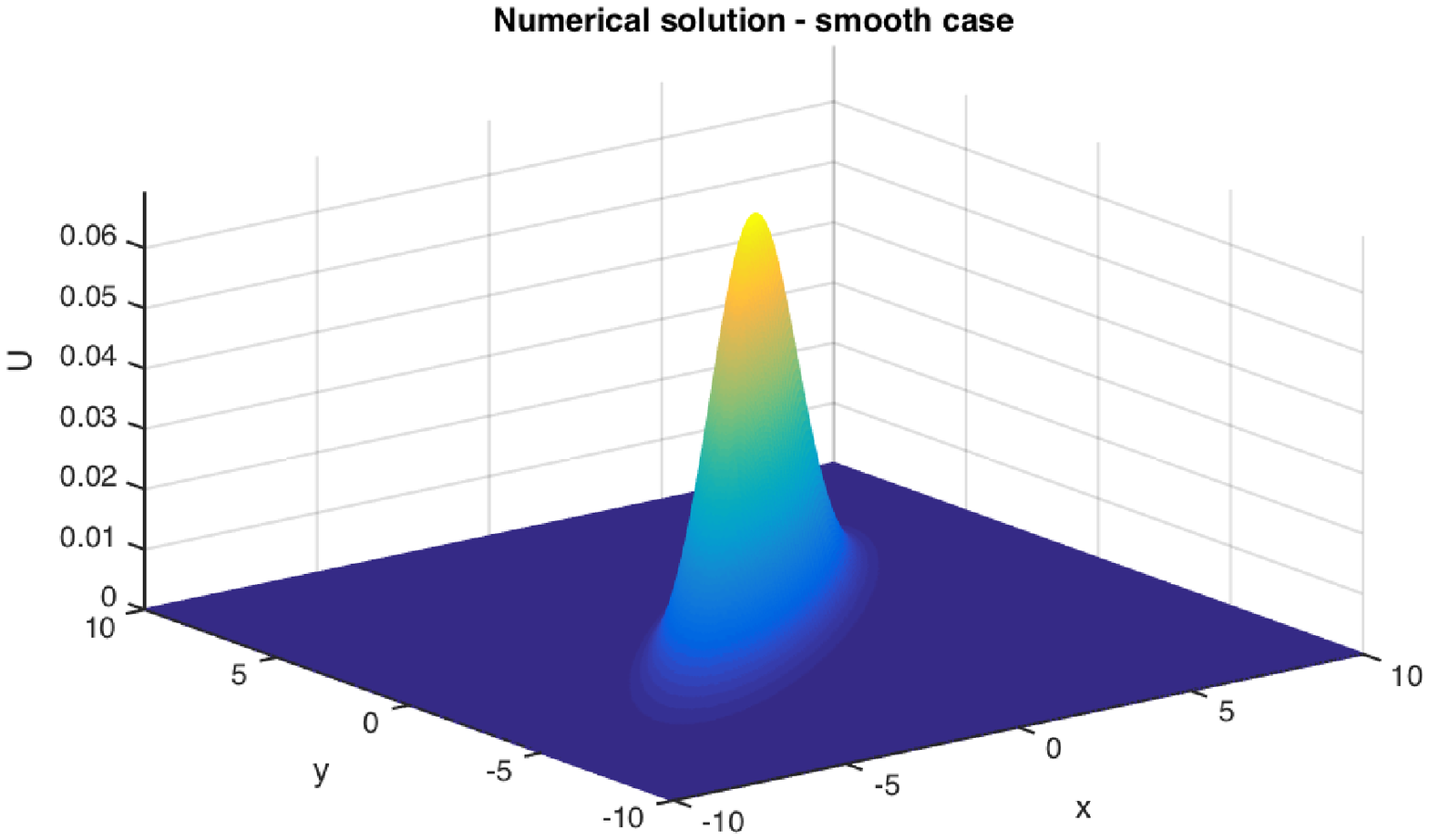}\hfill
\includegraphics[trim={0.5cm 1 1 1},clip, width=3.4in,height=2.2in]{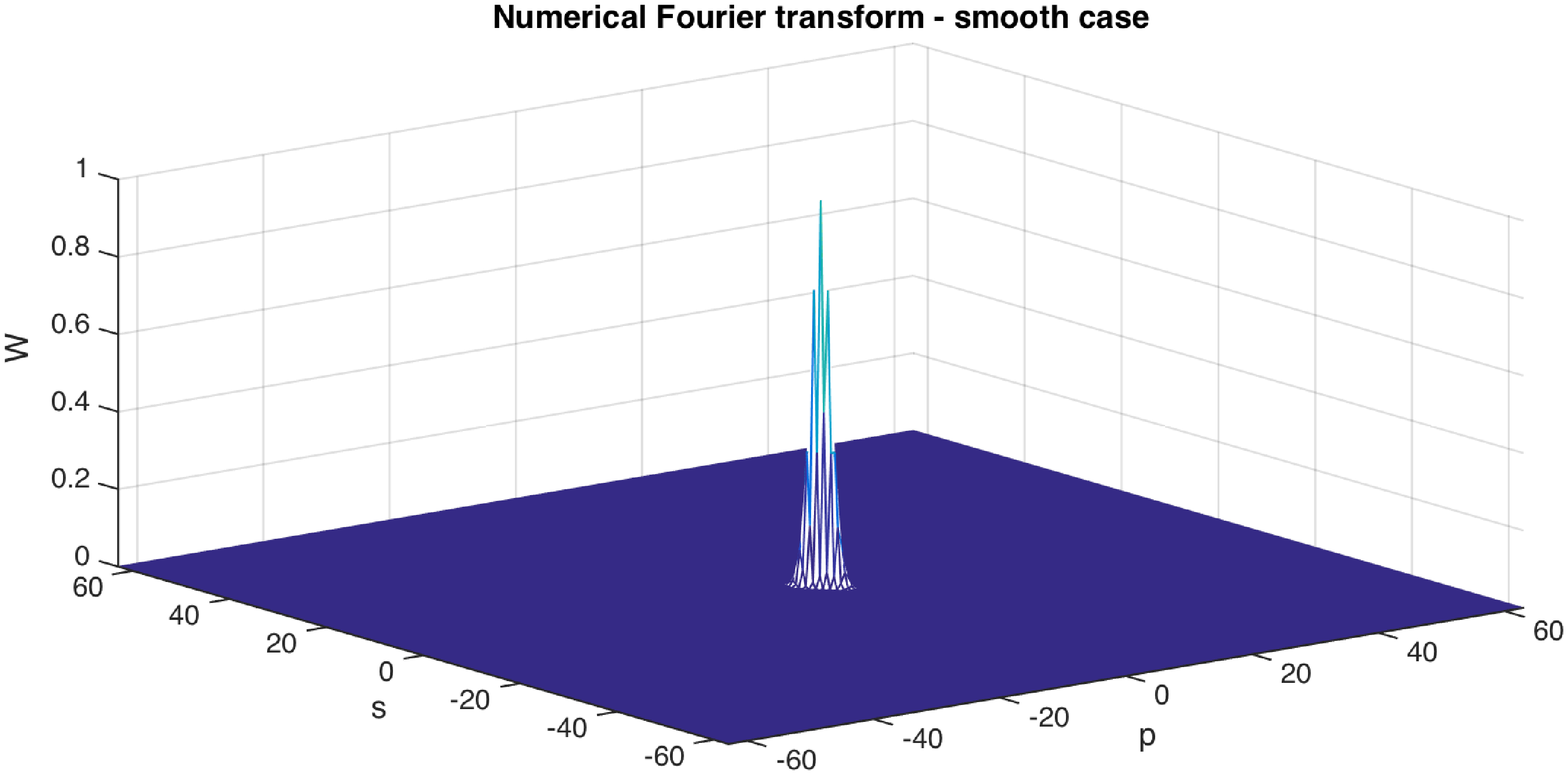}
\end{center}
\label{fig:smoothcentral}
\end{figure}

Extending the techniques from \cite{GC}, the time-discrete evolution of the discrete Fourier transform $W$ of the numerical solution $U$ is found to be governed by the recursion (see also Sections \ref{subsec:evol} and especially \ref{subsec:fullydiscr})
\begin{align*}
\left(1+4 \frac{\Delta t}{\Delta x^2} \sin^2\left(\frac{s \Delta x}{2}\right)\right) W^{n+1}(s,p) &= W^n(s,p) + \frac{\Delta t}{\Delta y} \sin(p \Delta y)
\frac{\partial W^{n+1}}{\partial s}(s,p),\qquad n \geq 0,\\
W^0(s,p) &= 1,
\end{align*}
for all $(s,p) \in [-\frac{\pi}{\Delta x}, \frac{\pi}{\Delta x}] \times [-\frac{\pi}{\Delta  y}, \frac{\pi}{\Delta y}]$. Setting $s=0$ and $p=\pm \pi/\Delta y$ we deduce that
\begin{eqnarray*}
W^n(0,\pm \pi/\Delta y) = 1 \qquad \forall\, n\ge 1.
\end{eqnarray*}
This finding is in line with our numerical simulations, and suggests that the numerical solution will not converge to the analytical solution
as $\Delta x, \Delta y \rightarrow 0$. To reconcile this evidence with
the results in \cite{PORRETTA}, we compute a numerical approximation to the solution at $T=2$, starting with the exact (smooth) solution (see (\ref{toyexact})) as initial datum at $t=1$.
The numerical solution and its discrete Fourier transform are shown in Fig.~\ref{fig:smoothcentral}.

As the initial datum in this case does not have high wavenumber components, the numerical solution approximates the analytical solution well. Indeed, the maximum error is around
$1.8 \times 10^{-4}$ (solution $\approx 0.07$), compared to $0.28$ (solution $\approx 0.3$) for the Dirac case.

Replacing the central $y$-difference with an upwind $y$-difference is observed to produce a convergent sequence of numerical approximations to the analytical solution in the case of a Dirac initial datum, just as for smooth
initial data, but such a finite difference scheme is only of first-order accuracy with respect to $\Delta y$. We shall therefore propose in the next section a numerical scheme applied to the analytical $y$-Fourier transform, which does not suffer from this shortcoming.

The remainder of the article is structured as follows.
In Section \ref{sec:scheme} we introduce the semidiscrete Fourier scheme and formulate the discretization of the model problem (\ref{toypde}), (\ref{toydelta}); in Section \ref{sec:analysis} we prove that the scheme is second-order convergent in $\Delta x$ in the space-time $\ell^\infty$ norm.
Next, in Section \ref{sec:fourier-sol}, we describe aspects of the numerical implementation including
the exponential convergence in $\Delta p$, the mesh size of the Fourier variable in the $y$-direction,
and present numerical results for the model problem (\ref{toypde}), (\ref{toydelta}).
In Section \ref{sec:hedge}, we discuss the application of the method to a financial hedging problem, while Section \ref{sec:multi} is concerned
with the convergence analysis of the extension of the scheme to the case of $n=2$, i.e. the diffusion operator acts in two spatial directions,
$x=(x_1,x_2)$, while in the third spatial direction, $y$, there is only a transport term present in the equation; in other words, there is no diffusion in the $y$-direction. Our findings are summarized in the concluding section.

\section{The semidiscrete Fourier scheme}
\label{sec:scheme}


The Fourier transform in the $y$-direction (briefly, $y$-FT) is defined by
\[v(x,p,t) := \int^{\infty}_{-\infty}u(x,y,t)\,{\rm e}^{i py}\dd y, \qquad x \in \mathbb{R},\; p \in \mathbb{R}, \; t> 0.\]
where it is tacitly understood that the function $y \in \mathbb{R} \mapsto u(x,y,t) \in \mathbb{R}$
is an element of $L^1(\mathbb{R})$ for all $(x,t) \in \mathbb{R}^n\times (0,\infty)$. The Fourier transform in the $y$-direction
of the initial Dirac measure $\delta(x-x_0) \otimes \delta(y-y_0)$ at $t=0$ is to be understood in the sense of distributions, and is equal to $\delta(x-x_0) \otimes {\rm e}^{ipy_0}$.

Application of the Fourier transform to \eqref{generalpde} yields
\[   \frac{\partial v}{\partial t}  -i pc(x,t) v =\mathcal{L}    v  \qquad x \in \mathbb{R},\; p \in \mathbb{R}, \; t> 0,\]
which is a family of PDEs in $x$ and $t$, parametrized by $p$, for the $y$-FT: $v(x,p,t)$. We then discretize the operator $\mathcal{L}$ from (\ref{Lop})
in the $x$-direction(s) by means of a standard finite difference scheme, using equally spaced grid points, with spacing $\Delta x$, but we keep the time variable continuous
for the moment at least. This leads to a system of ordinary differential equations (ODEs) indexed by the $x$-grid points, $x_j$, and also by the Fourier
wavenumber in the $y$-direction, which we are denoting by $p$, for the function $V_j(p,t)$. After solving this system of ODEs (in practice numerically, for a finite set of grid points $x_j$ subject to an artificial/numerical Dirichlet boundary condition at the `far-field', and a finite set of $y$-wave numbers, $\{p_l\}_{l=l_{min}}^{l_{max}})$, we use the inverse $y$-FT to obtain an approximate solution, $U(x_j,y,t)$, in the original co-ordinates as
\begin{eqnarray}
\label{inv-trunc}
U(x_j,y,t) =  \frac{1}{2\pi} \int^{2 t/(\Delta x)^r}_{-2 t/(\Delta x)^r}
V_j(p,t)\, {\rm e}^{-i y p} \dd p.
\end{eqnarray}
%
The truncation of the $p$-integration range in the inversion (\ref{inv-trunc}) from $p \in (-\infty,\infty)$ to the truncated range $p \in (-2 t/(\Delta x)^r, 2 t/(\Delta x)^r)$, for
a suitable $r>0$, is
dictated by the practical requirement to carry out numerical integration over a finite range.
The scaling with $2t$ simplifies the numerical analysis, but of course in practice any suitable scaling would be chosen empirically.


To find the approximation $U_{j,k}(t)$ to $u(x_j,y_k,t)$ numerically, we apply a uniformly spaced and equally weighted quadrature rule to (\ref{inv-trunc}) and obtain
\begin{eqnarray}
\label{numinvft}
U_{j,k}(t) := \frac{\Delta p}{2\pi}    \sum^{l_{max}}_{l=l_{min}}   V_{j}(p_l,t)\,{\rm e}^{-i y_k p_l}
\end{eqnarray}

\noindent for a given $k$, $l_{max}>0$, $l_{min}<0$, 
$l_{max} \Delta p = - l_{min} \Delta p = 2 t \Delta x^{-r}$. We will also denote $n_p := l_{max}-l_{min}+1$.
The numerical results will be seen to exhibit exponential convergence of the $p$-quadrature (see below and \cite{EXTTRAP}).
For an efficient implementation of (\ref{numinvft}), if the solution is needed for several values of $k$, the Fast Fourier Transform (FFT) can be used.

\medskip

We illustrate the method by applying it to the Cauchy problem
\begin{alignat}{2}
\label{toy2}
\frac{\partial u}{\partial t}+x\frac{\partial u}{\partial y}&=\frac{\partial^2 u}{\partial x^2}, &&\qquad (x,y,t) \in \mathbb{R}\times\mathbb{R}\times (0,T], \\
u(x,y,0) &= \delta (x)\otimes \delta (y), &&\qquad (x,y) \in \mathbb{R} \times\mathbb{R},
\label{toyic2}
\end{alignat}
with the aim to analyze the stability and accuracy of the numerical scheme we develop below.

Applying the $y$-FT to (\ref{toy2}) we get
%
\begin{alignat}{2}
\label{toyftpde}
\frac{\partial v}{\partial t}-i pxv &=\frac{\partial^2 v}{\partial x^2},
&&\qquad (x,p,t) \in \mathbb{R}\times\mathbb{R}\times (0,T], \\
u(x,p,0) &= \delta (x), &&\qquad (x,p) \in \mathbb{R} \times\mathbb{R},
\label{toyicft}
\end{alignat}
%
and we then discretize this one-parameter family of Cauchy problems in the $x$-direction only, using central differencing with spacing $\Delta x>0$, to obtain, for $x_j = j \Delta x$, $j \in \mathbb{Z}$,
\begin{align} \label{eq:SEMIDSCHEME}
\frac{\partial V_{j}(p,t)}{\partial t}-i p x_jV_{j}(p,t)=\frac{V_{j+1}(p,t)-2V_{j}(p,t)+V_{j-1}(p,t)}{\Delta x ^2},  \qquad p \in \mathbb{R}, \; t \in (0,T],
\end{align}
so that the function $V_{j}(p,t)$, which approximates $v(x_j,p,t)$, satisfies this equation for each $j\in \mathbb{Z}$ and for all $(p,t)
\in \mathbb{R}\times (0,T]$. The initial condition we use for this parametrized ODE system is
\begin{eqnarray}
\label{ic-semift}
V_{j}(p,0) = \left\{
\begin{array}{cl}
0 & \text{ for } j \ne 0, \\
\frac{1}{\Delta x} & \text{ for } j = 0,
\end{array}\right.
\end{eqnarray}
for all $p \in \mathbb{R}$, which approximates (\ref{toyicft}).

\section{Analysis of the numerical method for the stylized problem}
\label{sec:analysis}

In this section, we will prove the following theorem, which is one of our main results.
\begin{theorem}
\label{mainthm}
Let $u$ be the solution to the Cauchy problem \eqref{toy2}, \eqref{toyic2}, and let $U$ be given by \eqref{inv-trunc}, \eqref{eq:SEMIDSCHEME}, \eqref{ic-semift}. Then, for any $r>0$ in \eqref{inv-trunc},
\[
U(x_j,y,t) - u(x_j,y,t) \;= \; C(x_j,y,t) \, \Delta x^2 + o(\Delta x^2),\qquad j \in \mathbb{Z},\; y \in \mathbb{R},\; t>0,
\]
where
\[
C(x_j,y,t) := 
 \left[\frac{t}{2 \pi^2} \left( \frac{1}{4!} \left(\frac{\partial}{\partial x}+ \frac{t}{2} \frac{\partial}{\partial y} \right)^4 +
\frac{1}{5!} \left(\frac{t}{2} \frac{\partial}{\partial y} \right)^4 \right) u\right]\Bigg|_{(x_j,y,t)}.
\]
\end{theorem}

\subsection{The time evolution of the numerical double transform}
\label{subsec:evol}

To investigate the stability and accuracy of the numerical scheme, we use techniques motivated by those in \cite{GC} and \cite{CRAW}.
Thus, given a set of values $\{f_j\}_{j \in \mathbb{Z}} \in \ell^1(\mathbb{Z})$, on a uniformly-spaced grid $\{x_j\}_{j \in \mathbb{Z}}$ of spacing $\Delta x>0$
on $\mathbb{R}$, we consider the (semi-)discrete Fourier transform
\[
\widehat{f}(s) := \Delta x\sum^{\infty}_{j=-\infty}  f_j\,{\rm e}^{i sx_j},\qquad s \in \bigg[-\frac{\pi}{\Delta x}, \frac{\pi}{\Delta x}\bigg],
\]
and the inverse of this transform,
\[
f_j=\frac{1}{2\pi}\int^{\frac{\pi}{\Delta x}}_{-\frac{\pi}{\Delta x}}  \widehat{f}(s)\,{\rm e}^{i sx_j}   \dd s, \qquad j \in \mathbb{Z},
\]
(see, e.g.\ \cite{larsson2008partial}).
We note that the method of analysis here is specific to equation (\ref{toy2}) and its higher-dimensional variants (see Section \ref{sec:multi}),
but the numerical algorithm itself is not, as we demonstrate in the numerical example of Section \ref{sec:hedge}.

Then, by applying the semidiscrete $x$-FT to the system of ODEs (\ref{eq:SEMIDSCHEME}), we find
\begin{align*}
\frac{\partial W(s,p,t) }{\partial t} -p\frac{\partial W(s,p,t) }{\partial s}
&= \sum^{\infty}_{j=-\infty} \frac{\partial V_{j}(p,t)}{\partial t}\,  {\rm e}^{i sx_j} \Delta x - i p \sum^{\infty}_{j=-\infty}x_jV_{j}(p,t)\,{\rm e}^{i sx_j} \Delta x \\
&= \sum_{j=-\infty}^{\infty}\frac{V_{j+1}(p,t)\,{\rm e}^{i sx_j}-2V_{j}(p,t)\,{\rm e}^{i sx_j}+V_{j-1}(p,t)\,{\rm e}^{i sx_j}}{\Delta x ^2}\,  \Delta x  \\
&=W(s,p,t)\,\frac{4}{\Delta x ^2}\sin^2{\bigg(\frac{s\Delta x}{2}\bigg)},\qquad (s,p) \in \bigg[-\frac{\pi}{\Delta x}, \frac{\pi}{\Delta x}\bigg] \times \mathbb{R},\; t>0.
\end{align*}
As $W(s,p,0)=1$, it follows by the method of characteristics that $W(s,p,t)>0$ for all $t>0$, and therefore
\[
\frac{\partial \log{W} }{\partial t} -p\frac{\partial\log{W}  }{\partial s}  =\frac{4}{\Delta x ^2}\sin^2{\bigg(\frac{s\Delta x}{2}\bigg)},
\qquad (s,p) \in \bigg[-\frac{\pi}{\Delta x}, \frac{\pi}{\Delta x}\bigg] \times \mathbb{R},\; t>0.
\]
%

%
%
%
%

In contrast, taking the $x$-FT of (\ref{toyftpde}), the true double Fourier transform satisfies
\begin{equation*} \label{eq:w}
 \frac{\partial w }{\partial t}-p\frac{\partial w}{\partial s}=-s^2w, \qquad (s,p) \in \mathbb{R}^2,\; t>0.
\end{equation*}
Since $w(s,p,0)=1$, and therefore, by the method of characteristics $w(s,p,t)>0$ for all $t>0$, we have that
\begin{equation} \label{eq:logw}
   \begin{aligned}
     \frac{\partial \log{w} }{\partial t}-p\frac{\partial \log{w}}{\partial s}=-s^2, \qquad (s,p) \in \mathbb{R}^2,\; t>0.
   \end{aligned}
\end{equation}
Then, by defining
\[
Z(s,p,t)=\log{\left( \frac{W(s,p,t)}{w(s,p,t)} \right )},\qquad (s,p) \in \bigg[-\frac{\pi}{\Delta x}, \frac{\pi}{\Delta x}\bigg] \times \mathbb{R},\; t\geq 0,
\]
we find that
\[
   \frac{\partial Z(s,p,t) }  {\partial t} -p\frac{\partial Z(s,p,t)  }    {\partial s}  =  s^2-  \frac{4}{\Delta x ^2}\sin^2{\bigg(\frac{s\Delta x}{2}\bigg)}, \qquad (s,p) \in \bigg[-\frac{\pi}{\Delta x}, \frac{\pi}{\Delta x}\bigg] \times \mathbb{R},\; t>0,
\]
where
\[
g(s) =  s^2-  \frac{4}{\Delta x ^2}\sin^2{\bigg(\frac{s\Delta x}{2}\bigg)}, \qquad s \in \bigg[-\frac{\pi}{\Delta x}, \frac{\pi}{\Delta x}\bigg],
\]
and $Z(s,p,0)=0$. We can solve for $Z(s,p,t)$ to obtain
\[
Z(s,p,t)= \frac{1}{p} \int^{s+pt}_{s}  g(\sigma) \dd\sigma, \qquad (s,p) \in \bigg[-\frac{\pi}{\Delta x}, \frac{\pi}{\Delta x}\bigg] \times \mathbb{R},\; t \geq 0.
\]
Finally, we have
\begin{equation}
\label{eq:Ww}
W(s,p,t) =w(s,p,t) \exp{\left(\frac{1}{p} \int^{s+pt}_{s}  g(\sigma) \dd\sigma\right  )},\qquad (s,p) \in \bigg[-\frac{\pi}{\Delta x}, \frac{\pi}{\Delta x}\bigg] \times \mathbb{R},\; t\geq 0,
\end{equation}
an expression for the numerical double transform, $W$, in terms of the true double transform, $w$, with the exponential factor
on the right-hand side of \eqref{eq:Ww} reflecting the error introduced by the finite difference approximation in the $x$-direction.
In fact, one can solve (\ref{eq:logw}) with initial datum $w(s,p,0)=1$ to obtain
\begin{eqnarray}
\label{toydoubleft}
w(s,p,t)=\exp{\left(-s^2t-spt^2-\frac{1}{3}p^2t^3\right)}, \qquad (s,p) \in \mathbb{R}^2,\; t \geq 0.
\end{eqnarray}

A key observation is that it is more convenient to restate the solution in terms of the variables
\begin{eqnarray*}
\eta \;=\; s+\frac{pt}{2},\qquad
\xi \;=\; \frac{pt}{2},
\end{eqnarray*}
in Fourier space, a manifestation of the mixing of Fourier modes due to the lack of commutativity of the differential operators in $x$ and $y$ in (\ref{toy2}).
We will refer to these variables in Fourier space as wave numbers, and
it is only in these new variables that suitable wave number regimes can be defined. Indeed, in these new variables, we get
\begin{align}
\nonumber
W(s,p,t)
&=w(s,p,t) \exp{\left(\frac{1}{p} \int^{s+pt}_{s}  \sigma^2 \dd\sigma\right  )}  \exp{\left(-\frac{1}{p} \int^{s+pt}_{s} \frac{4}{\Delta x^2} \sin^2{\bigg(\frac{\sigma\Delta x}{2}\bigg)} \dd\sigma\right  )} \notag\\
\nonumber
&=w(s,p,t) \,w(s,p,t)^{-1} \exp{\left(-\frac{2}{p\Delta x^2} \int^{s+pt}_{s}(1- \cos{(\sigma\Delta x)}) \dd\sigma\right  )}\notag \\
&=\exp{ \bigg(-\frac{ 2t}{\Delta x^2}(1-\sinc{(\xi\Delta x)}\cos{(\eta\Delta x)}) \bigg)},
\qquad (s,p) \in \bigg[-\frac{\pi}{\Delta x}, \frac{\pi}{\Delta x}\bigg] \times \mathbb{R},\; t \geq 0,
\label{numdouble}
\end{align}
where sinc is the unnormalized sinc function, defined as follows (see, for example, \cite{OLVER}):
\[
\sinc x    := \frac{\sin x}{x}  \qquad \mbox{for $x \in \mathbb{R}\setminus\{0\}$ and $\sinc 0 :=1$}.
\]
We also have
\[
w(\eta,\xi,t)=\exp{\Big(-\eta^2t-\frac{1}{3}\xi^2t\Big)}, \qquad (\eta,\xi) \in \mathbb{R}^2,\; t \geq 0,
\]
in the new variables, $\eta$ and $\xi$.

\subsection{The wave number regimes} \label{sect:SEMIDANAL}
\begin{figure}
\begin{center}
\setlength{\unitlength}{0.5cm}
\begin{picture}(20,15)(0,0)
\put(15,7.5){\vector(1,0){4}}
\put(19.75,7.4){$\eta$}
\put(10,11.5){\vector(0,1){3}}
\put(10,15){\makebox(0,0){$\xi$}}
\linethickness{.4mm}
\put(9,11.5){\line(1,0){10}}
\put(6,8.5){\line(1,0){10}}
\put(4,6.5){\line(1,0){10}}
\put(1,3.5){\line(1,0){10}}
\put(9,6.5){\makebox(2,2){$\Omega_1$}}
\put(9,6.5){\line(0,1){2}}
\put(11,6.5){\line(0,1){2}}
\put(6,6.5){\makebox(2,2){$\Omega_2$}}
\put(12,6.5){\makebox(2,2){$\Omega_2$}}
\put(11.5,9){\makebox(2,2){$\Omega_3$}}
\put(6.5,4){\makebox(2,2){$\Omega_3$}}
\put(4,10){\makebox(2,2){$\Omega_4$}}
\put(15,3){\makebox(2,2){$\Omega_4$}}
\linethickness{1mm}
\put(1,3.5){\line(1,1){8}}
\put(19,11.5){\line(-1,-1){8}}
\end{picture}
\vspace{-2 cm}
\end{center}
\caption{Schematic representation of the wave number regimes defined by (\ref{omega1}) to (\ref{omega4}).}
\label{fig-regimes}
\end{figure}
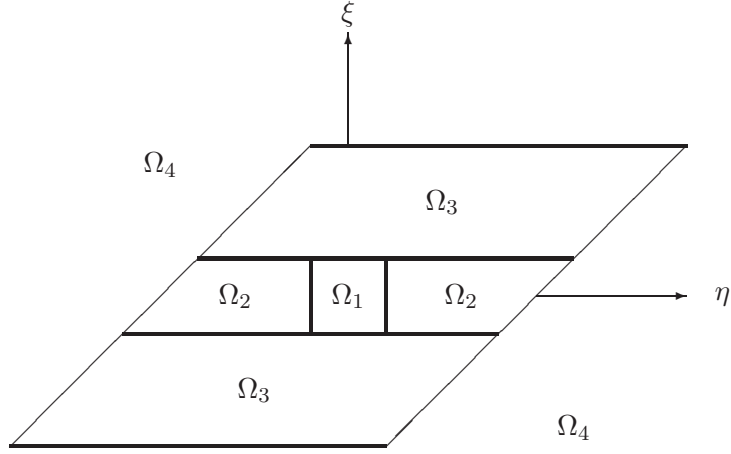

We decompose $\mathbb{R}^2$ into suitable wave number regimes in $\eta$ and $\xi$, for some $0<q<1$ and $r>0$,
\begin{align}
\label{omega1}
\Omega_1 &:= [ -\Delta x^{q-1},\Delta x^{q-1}   ] \times [   -\Delta x^{q-1},\Delta x^{q-1}    ],  \\
\Omega_2 &:=
\left\{ (\eta,\xi): -\Delta x^{q-1} \le \xi \le \Delta x^{q-1},\; -\frac{\pi}{\Delta x}+\xi \le \eta \le -\Delta x^{q-1}  \vee
\Delta x^{q-1} \le \eta \le \frac{\pi}{\Delta x}+\xi \right\}, \\
\Omega_3 &:=
\left\{ (\eta,\xi):
-\frac{\pi}{\Delta x}+\xi \le \eta \le  \frac{\pi}{\Delta x}+\xi,\;
-\Delta x^{-r} \le \xi \le -\Delta x^{q-1} \vee \Delta x^{q-1} \le \xi \le \Delta x^{-r}
\right\}, \\
\Omega_4 &:= \mathbb{R}^2 \backslash (\Omega_1\cup \Omega_2\cup \Omega_3),
\label{omega4}
\end{align}
which are also illustrated in Fig.~\ref{fig-regimes}.
By applying the inverse transforms to the error term we obtain an equation for the error in the $(x,y)$-plane for $t>0$, i.e.,
\begin{align*}
\setbraceratio{1}{3}
(2\pi)^2(U(x,y,t)-u(x,y,t))  &=  \int^{2 t \Delta x^{-r}}_{-2 t \Delta x^{-r}} \int^{\frac{\pi}{\Delta x}}_{-\frac{\pi}{\Delta x}}
(W(s,p,t)-w(s,p,t))\, {\rm e}^{-i sx}\,{\rm e}^{-i py} \dd s \dd p   \\
& \qquad - \int_{\Omega_4} w(s,p,t)\, {\rm e}^{-i sx}\, {\rm e}^{-i py} \dd s \dd p  \\
&=
\frac{2}{t}
\int^{\Delta x^{-r}}_{-\Delta x^{-r}} \int^{\frac{\pi}{\Delta x}+\xi}_{-\frac{\pi}{\Delta x}+\xi} (W(\eta,\xi,t)-w(\eta,\xi,t))\,
{\rm e}^{-i (\eta-\xi)x}\, {\rm e}^{-i \frac{2\xi y}{t}} \dd \eta \dd \xi
\\
& \qquad - \frac{2}{t} \int_{\Omega_4} w(\eta,\xi,t)\, {\rm e}^{-i(\eta-\xi)x}\, {\rm e}^{-i \frac{2\xi y}{t}} \dd\eta \dd\xi.
\end{align*}
Here $x$ is any $x$-grid node. We now define
\begin{align*}
I_1(x,y,t) &:=  \int_{\Omega_1} (W(\eta,\xi,t)-w(\eta,\xi,t))\,
{\rm e}^{-i(\eta-\xi)x}\, {\rm e}^{-i\frac{2\xi y}{t}} \dd\eta \dd\xi, \\
I_k(x,y,t) &:=  \int_{\Omega_k} W(\eta,\xi,t)\, {\rm e}^{-i(\eta-\xi)x}\, {\rm e}^{-i\frac{2\xi y}{t}} \dd\eta \dd\xi, \qquad k =2,3, \\
I_4(x,y,t) &:=  \int_{\mathbb{R}^2 \backslash \Omega_1} w(\eta,\xi,t)\, {\rm e}^{-i(\eta-\xi)x}\, {\rm e}^{-i\frac{2\xi y}{t}} \dd\eta \dd\xi.
\end{align*}
Except for the joint low wavenumber regime, we will perform separate calculations on $W$ and $w$.
We will find that all but the first term can be made exponentially small, while the first term is $O(\Delta x^2)$.
We begin by considering the joint low wavenumber analysis for $W$ and $w$.

\subsection{Joint low wavenumbers (region $\boldsymbol{\Omega}_{\mathbf{1}}$)}

We have to determine an exponent $q \in (0,1)$ such that, in $\Omega_1$,
\[
|\eta \Delta x|<\Delta x^q  \rightarrow 0 \qquad \text{ and } \qquad
|\xi \Delta x|<\Delta x^q   \rightarrow 0,
\]
and certain expansions can be usefully truncated. By straightforward Taylor expansion,
\begin{align*}
-\frac{2t}{\Delta x^2}(1-\sinc{(\xi\Delta x)}\cos{(\eta\Delta x)})
&= -\eta^2t-\frac{1}{3} \xi^2t + \frac{2t}{4!} \eta^4\Delta x^2 +\frac{2t}{5!} \xi^4\Delta x^2 \\
& \quad\, +\, O(\eta^6\Delta x^4)+O(\xi^6\Delta x^4).
\end{align*}
Our objective is to choose $q \in (0,1)$ so that the remainder terms, resulting from approximating the integrand of $I_1$ by its Taylor series expansion, are $o(\Delta x^2)$. To this end, we write
\begin{align}
\nonumber
I_1(x,y,t)
&= \int_{\Omega_1}
w\, \bigg (\exp{\left(\frac{2t}{4!} \eta^4\Delta x^2 +\frac{2t}{5!} \xi^4\Delta x^2+o(\Delta x^2)\right)} -1 \bigg)\, {\rm e}^{-i(\eta-\xi)x}
\,{\rm e}^{-i\frac{2\xi y}{t}} \dd \eta \dd \xi \\
\label{extension}
&=   \Delta x^2 \int_{\Omega_1}
w
\left (\frac{2t}{4!} \eta^4\Delta x^2 +\frac{2t}{5!} \xi^4\Delta x^2\right) {\rm e}^{-i(\eta-\xi)x}\, {\rm e}^{-i\frac{2 \xi y}{t}} \dd \eta \dd\xi +
o(\Delta x^2) \\
\nonumber
&=   \Delta x^2 \int_{\mathbb{R}^2}
w
\left (\frac{2t}{4!} \eta^4\Delta x^2 +\frac{2t}{5!} \xi^4\Delta x^2\right) {\rm e}^{-i(\eta-\xi)x} \,{\rm e}^{-i\frac{2 \xi y}{t}} \dd \eta \dd\xi
+o(\Delta x^2) \\
&=   \Delta x^2\, F(x,y,t) +o(\Delta x^2),
\nonumber
\end{align}
where
\begin{align*}
F(x,y,t) &:= \int_{\mathbb{R}^2}
w \left (\frac{2t}{4!} \eta^4 +\frac{2t}{5!} \xi^4\right)  \,{\rm e}^{-i(\eta-\xi)x} \,{\rm e}^{-i\frac{2 \xi y}{t}} \dd \eta \dd\xi \\
&= \frac{t}{2} \int_{\mathbb{R}^2}
w \left (\frac{2t}{4!} \left(s+\frac{p t}{2} \right)^4 +\frac{2t}{5!} \left(\frac{p t}{2} \right)^4\right)  \,{\rm e}^{-i s x} \,{\rm e}^{-i p y} \dd s \dd p \\
&= t^2 \left( \frac{1}{4!} \left(\frac{\partial}{\partial x}+ \frac{t}{2} \frac{\partial}{\partial y} \right)^4 +
\frac{1}{5!} \left(\frac{t}{2} \frac{\partial}{\partial y} \right)^4 \right) u.
\end{align*}
We are able to replace $\Omega_1$ by $\mathbb{R}^2$ in (\ref{extension})  to $o(\Delta x^2)$ because $w$ decays exponentially
in $\eta$ and $\xi$ (see also Section \ref{sec:loc}). The last step uses the relation between the Fourier transform of a smooth function and the transforms of its derivatives.
We also require
\[ \Delta x^2( \eta^8+\xi^8        )  \rightarrow 0   \]
for the remainders in the Taylor expansion of the exponential to be $o(\Delta x^2)$. We can therefore take $q \in (\frac{3}{4}, 1)$ to define the joint low wavenumber regime.

Having dealt with the region $\Omega_1$ in the $(s,p)$-plane, we move on to consider the remaining terms involving $W$. There are two cases to discuss, corresponding
to low $\xi$-wavenumbers and high $\eta$-wavenumbers (region $\Omega_2$), and to high $\xi$-wavenumbers (region $\Omega_3$).

\subsection{Low $\boldsymbol{\xi}$-wavenumbers and high $\boldsymbol{\eta}$-wavenumbers (region $\boldsymbol{\Omega}_{\mathbf{2}}$)}
\label{low-high}
For these wavenumbers we have
$\eta\Delta x \in [-\pi+\xi\Delta x,\pi+\xi\Delta x]$ and in this interval $\eta=0$ is the only solution to $\cos{(\eta\Delta x)}=1$
since $|\xi \Delta x| \le \Delta x^{q}$. The following inequality is valid for $\theta \in [-\pi, \pi]$:
\[
\sin^2{ \bigg(\frac{\theta}{2} \bigg)} \ge  \bigg(\frac{\theta}{\pi} \bigg)^2.
\]
However since we wish to take $\theta = \eta\Delta x$ with $\theta \in [-\pi+\xi\Delta x,\pi+\xi\Delta x]$ and $|\xi \Delta x| \le \Delta x^{q}$, we shall use instead
the weaker inequality
\[
\sin^2{ \bigg(\frac{\theta}{2} \bigg)} \ge  \frac{1}{4}\bigg(\frac{\theta}{\pi} \bigg)^2,
\]
which is valid for all such $\theta$, provided that $\Delta x$ is sufficiently small (whereby also $|\xi \Delta x| \le \Delta x^{q}$ is sufficiently small).
Hence,
\[
\cos{(\eta\Delta x)}  =  1- 2\sin^2{ \bigg(\frac{\eta\Delta x}{2} \bigg)} \le 1-\frac{1}{2} \bigg (\frac{\eta\Delta x }{\pi} \bigg)^2 = 1-\frac{1}{2\pi^2} (\eta\Delta x )^2.
\]
We can then choose $\alpha<\frac{1}{2\pi^2}$ such that for small enough $\Delta x$ we have
\[
\cos{(\eta\Delta x)} \leq 1-\alpha(\eta \Delta x)  ^2
\]
for $\eta\Delta x \in [-\pi + \xi \Delta x, \pi + \xi \Delta x]$, since $|\xi\Delta x|\leq \Delta x^{q} \rightarrow 0$ as $\Delta x\rightarrow 0$.
Also, for $\xi\Delta x$ small enough, we can ensure that
\[
\frac{1}{2}<\sinc{(\xi\Delta x)} \le 1.
\]
So then, for $\cos(\eta\Delta x) \geq 0$,
\begin{eqnarray}
\label{high-eta-est}
\alpha(\eta\Delta x)^2  \le  1-\sinc{(\xi\Delta x)} \cos{(\eta\Delta x)}.
\end{eqnarray}
If, on the other hand, $\cos(\eta \Delta x)\leq 0$, then, because we
always have $\alpha (\eta \Delta x)^2 \le \frac{1}{2\pi^2} (\pi + \Delta x^q)^2 \le 1$ for $\Delta x$ sufficiently small (more precisely, for
$\Delta x \leq [(\sqrt{2}-1)\pi]^{1/q}$), while the
right-hand side of \eqref{high-eta-est} is $\ge 1$ for $\cos(\eta \Delta x)\leq 0$, it once again follows that \eqref{high-eta-est} holds.
Hence,
\begin{eqnarray*}
 \exp{ \bigg(-\frac{2t}{\Delta x^2}(  1-\sinc{(\xi\Delta x)} \cos{(\eta\Delta x)}) \bigg)}   \le
 \exp{\bigg(-\frac{2t}{\Delta x^2}\alpha (\eta\Delta x)^2\bigg)}
 \le  \exp{(-2\alpha t\Delta x^{2(q-1)})}
\end{eqnarray*}
as $\eta\Delta x >\Delta x ^{q}$ and $\eta^2 >\Delta x ^{2(q-1)}$. Thus, the contribution to the integral satisfies
\begin{align*}
|I_2|&\le\int_{\Omega_2}W(\eta,\xi,t) \dd \eta \dd \xi 
= o(\Delta x^m)\qquad \mbox{as $\Delta x \rightarrow 0$}  \notag
\end{align*}
for all $m>0$. We therefore deduce that this contribution to the integral is exponentially small as $\Delta x \rightarrow 0$.

\subsection{High $\boldsymbol{\xi}$-wavenumbers (region $\boldsymbol{\Omega}_{\mathbf{3}}$)}
\label{sec:high}

\begin{figure}[H]
\caption{Plot of $\sinc x=\frac{\sin{x}}{x}$ and its relation to other functions needed in 
 Section \ref{sec:high}.}
\centering
\includegraphics[width=5.5in,height=2.5in]{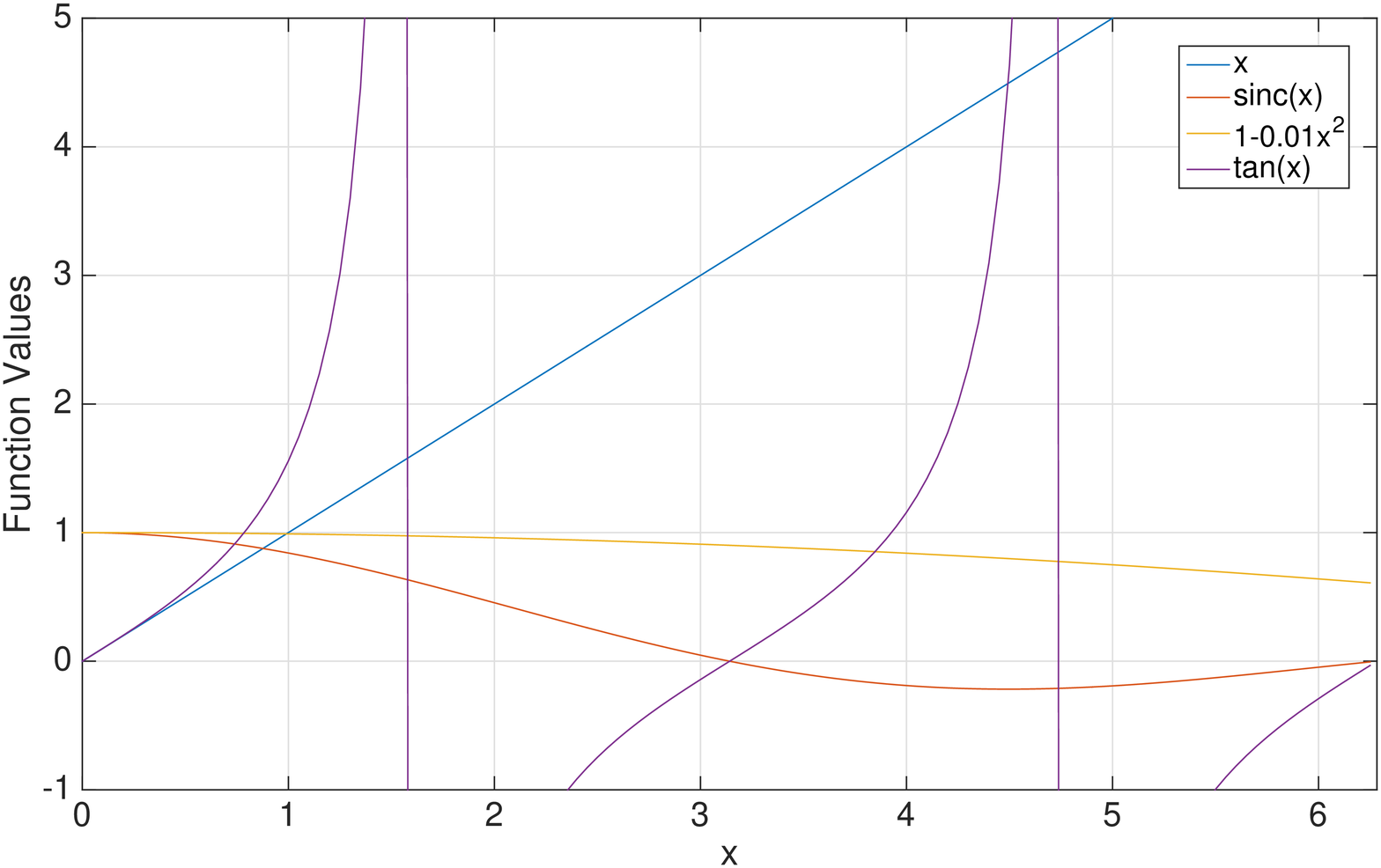}
\label{fig:SINCPLOT}
\end{figure}

We observe that $\xi \mapsto \sinc{(\xi\Delta x)}$ is a decreasing function
from $\xi=\Delta x ^{q-1}$ to some $\xi=\xi_0$, if we make $\Delta x$ small enough.
At that first local minimum, $\xi_0$, one then has $\tan{(\xi_0\Delta x)}=\xi_0\Delta x$ and then
$|\sinc{(\xi_0\Delta x)}|=\alpha$ with $\alpha<0.3$, say, and for values of $\xi>\xi_0$, we have $|\sinc{(\xi\Delta x)}|<\alpha$.
This is illustrated in Figure \ref{fig:SINCPLOT}, with the symbol {\footnotesize{\sf{x}}} in the figure signifying $\xi \Delta x$.
We see that for $0<\xi\le \xi_0$ we also have (see again Figure \ref{fig:SINCPLOT})
\[
 \sinc(\xi\Delta x) \cos(\eta \Delta x) \le |\sinc{(\xi\Delta x)}|  <1-\frac{1}{100} (\xi\Delta x)^2,
\]
while for $\xi \ge \xi_0$ we have
\[
 \sinc(\xi\Delta x) \cos(\eta \Delta x) \le | \sinc(\xi\Delta x) | \le \alpha\, ( < 0.3).
\]
So, for  $\xi\Delta x \ge  \Delta x^{q}$, the following inequality holds:
\[
1- \sinc(\xi\Delta x)\cos(\eta \Delta x) \ge \frac{1}{100}\min\left(\Delta x^{2q}, 100(1-\alpha) \right)\!.
\]
Therefore,
\begin{align*}
\exp{\bigg(-\frac{2t}{\Delta x^2}(  1-\sinc{(\xi\Delta x)} \cos{(\eta\Delta x)}        ) \bigg)}
&\le     \exp{ \bigg(  -\frac{t}{50\Delta x^2}\min{(\Delta x^{2q},100 (1-\alpha))}  \bigg )} \\
&\le     \exp{ \bigg(  -\frac{t}{50\Delta x^{2(1-q)}}  \bigg)}
\end{align*}
for small enough $\Delta x$. Hence, similarly as before (cf. Section \ref{low-high}),
\begin{align*}
| I_3| 
\le  &(  2\Delta x^{-r}     ) ( 2\pi\Delta x^{-1}  )  \exp{ \bigg(  -\frac{t}{50} \Delta x^{-2(1-q)} \bigg )}  =o(\Delta x ^m)
\qquad \mbox{as $\Delta x \rightarrow 0$}
\end{align*}
for all $m>0$.  So this contribution to the integral is also exponentially small as $\Delta x \rightarrow 0$.\\

\begin{remark}
We note that for the high $\xi$-wavenumber integration range to be nonempty, we need that
\[
\Delta x^{-(1-q)}<\xi<\Delta x^{-r} ,
\]
i.e.\ $r>1-q$. If $r\le 1-q$, there is only a low $\xi$-wavenumber regime, but this is irrelevant for the computations, and 
quadratic convergence is still guaranteed for any $r  \in (0,1)$ because the contributions from outside $(-\Delta x^{-r}, \Delta x^{-r})$
decay exponentially for all $r$ as $\Delta x \rightarrow 0$.
\end{remark}

\begin{remark}
For any fixed $\Delta x$, the numerical double Fourier transform $W$ in (\ref{numdouble}) approaches a constant value of $\exp(-2t/\Delta x^2)$ as
$\xi \rightarrow \pm \infty$ (equivalently, $p \rightarrow \pm \infty$) independent of $\eta$ (or $s$), hence
\[ \lim_{p \rightarrow \pm \infty} V_j(p,t) = \lim_{p \rightarrow \pm \infty} \frac{1}{2\pi} \int_{-\frac{\pi}{\Delta x}}^{\frac{\pi}{\Delta x}}
W(s,p,t)\,{\rm e}^{i s x_j} \dd s = \left\{ \begin{array}{cc} \frac{1}{\Delta x}\,\exp(-2t/\Delta x^2) & \mbox{if $j=0$},\\
0 & \mbox{if $j \neq 0$}. \end{array} \right.\]
This implies that for $j=0$, the integrand appearing in (\ref{inv-trunc}), as a function of $p$, does not belong to $L^1(\mathbb{R})$,
which is yet another reason why instead of using the actual inverse Fourier transform in (\ref{inv-trunc}) to define $U(x_j,y,t)$, with
$\mathbb{R}$ as integration range, we have integrated over the compact interval
$[-2t \Delta x^{-r}, 2t \Delta x^{-r}]$. However, because $V_0(p,t)$ rapidly decays to zero with $\Delta x \rightarrow 0$ for $t>0$ as  $p \rightarrow \pm \infty$, the sequence of integrals over $[-2t\Delta x^{-r},2t\Delta x^{-r}]$ converges to the true inverse Fourier transform, as $\Delta x \rightarrow 0$.
\end{remark}

\subsection{Localization of the exact Fourier transform (region $\boldsymbol{\Omega}_{\mathbf{4}}$)}
\label{sec:loc}

Finally, we need to estimate the error contribution for the remaining terms, which involve  $w$.
We get
\begin{eqnarray*}
|I_4|\le \int^{\infty}_{-\infty}  \int^{-\Delta x^{q-1}}_{-\infty} \!\! w \, \dd\eta \dd\xi
+  \int^{\infty}_{-\infty} \int_{\Delta x^{q-1}}^{\infty} \!\! w \, \dd\eta \dd\xi
+ \int^{\infty}_{\Delta x^{q-1}}\int^{\infty}_{-\infty} \!\! w \, \dd\eta \dd\xi
+  \int_{-\infty}^{-\Delta x^{q-1}} \!\!\! \int^{\infty}_{-\infty}  \!\! w \, \dd\eta \dd\xi.
\end{eqnarray*}
As we have
\[
\int^{\infty}_{-\infty} w(\eta,\xi,t)\dd\eta  =\exp{ \bigg(-\frac{1}{3} \xi^2 \bigg)} \int^{\infty}_{-\infty} \exp{(-
\eta^2)} \dd\eta  = \sqrt{\pi} \exp{ \bigg(-\frac{1}{3} \xi^2 \bigg)}
\]
and
\[  \int^{\infty}_{-\infty} w(\eta,\xi,t)\dd\xi  =\exp{(-
\eta^2)} \int^{\infty}_{-\infty} \exp{ \bigg(-\frac{1}{3}
\xi^2 \bigg)} \dd\xi  = \sqrt{3\pi} \exp{(-
\eta^2)},
\]
it follows that
\begin{align*}
|I_4| &\le  \sqrt{3\pi}  \int^{-\Delta x^{q-1}}_{-\infty} \exp{(-
\eta^2)} \dd\eta    +   \sqrt{3\pi}\int_{\Delta x^{q-1}}^{\infty} \exp{(-
\eta^2)} \dd\eta       \\
& \quad\, + \sqrt{\pi}  \int^{\infty}_{\Delta x^{q-1}} \exp{ \bigg(-\frac{1}{3}
\xi^2 \bigg)} \dd\xi   + \sqrt{\pi} \int_{-\infty}^{-\Delta x^{q-1}} \exp{ \bigg(-\frac{1}{3}
\xi^2 \bigg)}  \dd\xi .
\end{align*}
Lemma 3 in the Appendix of \cite{GILES} implies that each of these integrals is $o(\Delta x^m)$,
for any $m>0$, and so $|I_4|$ is exponentially small as $\Delta x \rightarrow 0$.

Collecting the above results we find that
\begin{align*}
U(x,y,t)-u(x,y,t) 
&= \frac{1}{2 t \pi^2} \left( I_1(x,y,t) + I_2 + I_3 - I_4 \right) \\
&= \frac{1}{2 t \pi^2} \left( I_1(x,y,t) + o(\Delta x^m) \right) \\
&= \Delta x^2 \left[\frac{t}{2 \pi^2} \left( \frac{1}{4!} \left(\frac{\partial}{\partial x}+ \frac{t}{2} \frac{\partial}{\partial y} \right)^4 +
\frac{1}{5!} \left(\frac{t}{2} \frac{\partial}{\partial y} \right)^4 \right) u\right]\Bigg|_{(x,y,t)} + o(\Delta x^2) \quad \mbox{as $\Delta x \rightarrow 0$},
\end{align*}
which completes the proof of Theorem \ref{mainthm}.

\section{Computation of the semidiscrete Fourier solution}
\label{sec:fourier-sol}

In this section we present the results of applying the semidiscrete Fourier scheme to the toy model \eqref{toypde}, \eqref{toydelta}.
We compute the solution to  (\ref{eq:SEMIDSCHEME}), (\ref{ic-semift}) by using the matrix exponential
and we then use the trapezium rule
and the Fast Fourier Transform (FFT) to compute the values of $U(x_j,y_k,T)$ from (\ref{inv-trunc}).

\subsection{Solving the ODEs}

The ODE system (\ref{eq:SEMIDSCHEME}) can be written as
\begin{align} \label{eq:SEMID}
 \frac{\rm{d}}{{\rm d}t} {V}(p,t)  =M {V}(p,t),
\end{align}
where $V=(V_j)_{j}$, $M=M_1+ipM_2$ is a bi-infinite matrix, with
\[M_1=\frac{1}{\Delta x^2}
{
\arraycolsep=1.6pt\def\arraystretch{1}
\left(\begin{matrix}
\ldots & \ldots & \ldots & \ldots & \ldots & \ldots & \ldots &\\

 \ldots & -2 & ~~1  & ~~0 & ~~0 & ~~0 & \ldots\\
 \ldots &  ~~1 & -2& ~~1 & ~~0 & ~~0 & \ldots\\

 \ldots &  ~~0 & ~~1 & -2& ~~1 & ~~0 & \ldots\\
  \ldots &  ~~0 & ~~0 &~~1  & -2 & ~~1 & \ldots\\
  \ldots &   ~~0 & ~~0 & ~~0 & ~~1  & -2 & \ldots \\
  \ldots & \ldots & \ldots & \ldots & \ldots & \ldots & \ldots & \\
\end{matrix}\right)
}, \qquad
%
%
M_2=\Delta x
{
\arraycolsep=1.6pt\def\arraystretch{1}
\left(\begin{matrix}
\ldots & \ldots & \ldots & \ldots & \ldots & \ldots & \ldots &\\
\ldots &  -2& ~~0   & ~~0 & ~~0 & ~~0 & \ldots\\
\ldots &   ~~0  & -1 & ~~0 & ~~0 & ~~0 & \ldots\\

\ldots &  ~~0 & ~~0  & ~~0 & ~~0& ~~0 & \ldots\\
\ldots &   ~~0 & ~~0 & ~~0  & ~~1 & ~~0 & \ldots\\
\ldots &    ~~0 & ~~0 & ~~0 & ~~0  & ~~2 & \ldots \\
\ldots & \ldots & \ldots & \ldots & \ldots & \ldots & \ldots & \\
\end{matrix}\right).
}
\]
In our implementation, the problem is considered over a sufficiently large square domain $x \in (-L, L)$ with zero Dirichlet boundary
condition at $x=\pm L$; in our numerical experiment below we took $L$ large enough ($L=10$ or $L=20$) to ensure that the Dirichlet boundary condition has negligible influence on the values of the numerical solution at the final time of interest, $T>0$. Hence, the bi-infinite matrices $M_1$, $M_2$ and
$M$ are truncated to square matrices $\widetilde{M_1}$, $\widetilde{M}_2$ and $\widetilde{M}$ of a certain finite size (depending on the choice of $L$ and $\Delta x$). Since the matrix $\widetilde{M}$ is independent of $t$, the truncated counterpart of (\ref{eq:SEMID}) has the obvious solution:
\begin{align}
\label{matexp}
\widetilde{V}(p,t)= {V}(p,0)\,{\rm e}^{\widetilde{M}t},\qquad p \in \mathbb{R}, \; t>0,
\end{align}
where $V_j(p,0)=0$ for $j \ne 0$ and $V_0(p,0)=\frac{1}{\Delta x}$.

We use the Matlab {\tt expm} function for the matrix exponential, which is based on the scaling and squaring method (cf. \cite{higham2005scaling}).
To improve this part of the algorithm, one could exploit that the eigenvalues of the symmetric matrix $\Delta x^2 \widetilde{M}$ are contained in the Gershgorin discs with centers $-2 + i j p \Delta x^3$, and radii $\leq 2$, where $|j p \Delta x^3| \le x_{\max} \Delta x^{2-r}$.
We do not pursue this further here.

The exact solution (\ref{matexp}) in the form of a matrix exponential is only available because the coefficients of \eqref{toypde}
do not depend on time.
We present results for a time stepping method in Section \ref{subsec:fullydiscr}.


The solution of (\ref{eq:SEMID}) can, in principle, be found for any value of $p$ but, in practice we will only be able to calculate a finite number of values, $\widetilde{V}_{j}(p_k,T)$.
To compute the solution in the original variables, we numerically invert the $y$-FT, using the grid values, $\widetilde{V}_{j}(p_k,T)$. By using this method, we have avoided discretizing the $y$-partial derivative: in effect, this is replaced by a discretization in the $p$-direction.

\subsection{Convergence of the $\mathbf{p}$-discretization}

When we solve the Kolmogorov forward equation (KFE) using the semidiscrete Fourier method, we need to consider the effect of the numerical
inversion of the $y$-FT in the final step of the solution procedure. As the inversion is only approximate, we have to decide how to set the
parameters for the inversion, and, more specifically, how to choose the range of $p$-wave numbers and the $p$-step.
%

According to the
Euler--Maclaurin expansion of the error of the composite trapezium rule applied to a sufficiently smooth function $f$ over an interval $[a,b]$
(see, for example, \cite{suli}, pp.\ 213, Theorem 7.4),
\begin{eqnarray}\label{eml}
\sum_{j=1}^k c_j h^{2j} \left[f^{(2k-1)}(b) - f^{(2k-1)}(a) \right] - \left(\frac{h}{2}\right)^{2k} \int_a^b q_{2k}(\tau(p)) f^{(2k)}(p) \ {\rm d}p,
\end{eqnarray}
where $h = (b-a)/m$ is the uniform mesh size, $\tau$ is a piecewise linear `saw tooth' function with values in $[-1,1]$, and $q_{2k}$ and $c_j$, $j=1,\dots,k$, are computable polynomials and constants, respectively.

The integrand in \eqref{inv-trunc}, up to a $p$-integrable $O(h^2)$ term, is, for $x$, $y$ and $t$ fixed,
\begin{equation}
\label{toysingleft}
f(p)=V(x,p,t) \exp(-ipy) = \frac{1}{2\sqrt{\pi t}} \exp\left(- \frac{t^3}{12} p^2 + i \frac{xt}{2} p - \frac{x^2}{4t}  \right) \exp(-ipy),
\end{equation}
as one finds by taking the inverse Fourier transform (in $p$) of (\ref{toydoubleft}).

%

In our case, the integrand (\ref{toysingleft}) and its derivatives (in $p$) vanish so rapidly for large $p$ that by taking the integration limits to be
$a = -2t\Delta x^{-r}$ and $b = 2t\Delta x^{-r}$ the first term in \eqref{eml} becomes exponentially small in $\Delta x$ for any fixed $k$.
Then, taking $m=n_p:=l_{max}-l_{min} + 1$ in \eqref{numinvft} such that $h = \Delta p = 4t \Delta x^{-r}/n_p \rightarrow 0$, say $h=\Delta x^{r}$, the second term is $O(h^{2k})$ for any $k$.
Hence, the total error can be made exponentially small in $\Delta x$.
See \cite{EXTTRAP} for a detailed discussion of the behaviour of the trapezium rule for analytic integrands in a neighbourhood of the real line.

\subsection{Using the inverse FFT}

When inverting the $y$-FT, we use the formula (\ref{numinvft}) and we can take advantage of the low computational complexity
of the Inverse Fast Fourier Transform (IFFT). Specifically, we use the MATLAB procedure {\tt ifft}. In the MATLAB documentation \cite{MATLAB}, the IFFT
is defined by the discrete inverse Fourier transform as
\[
X_j=\frac{1}{N}\sum^N_{l=1}Y_l \omega_N^{-(j-1)(k-1)},
\]
where
\[
\omega_N= {\rm e}^{-2i\pi/N},
\]
and where both $X$ and $Y$ have length $N$. We have to relate our inversion formula to the definition of the inversion formula in MATLAB;
the details of this are described in the Appendix \ref{app:ifft}.

\subsection{Numerical experiments}

By taking the inverse Fourier transform (in $s$) of $V$ in (\ref{toysingleft}), one finds the exact solution to (\ref{toypde}),
\begin{eqnarray}
\label{toyexact}
u(x,y,t) = \frac{\sqrt{3} }{2\pi t^2  } \exp{\left(-\frac{1}{4t} x^2 \right )    }         \exp{\left(-\frac{3   }{ t^3 } \left(y-\frac{t}{2}  x\right)^2\right)}, \qquad (x,y) \in \mathbb{R}^2,\; t>0,
\end{eqnarray}
(see also \cite{KOLM34}), with which the numerical solution can be compared.

\begin{figure}[H]
\caption{Toy model: maximum error in the semi-discrete Fourier method against the number of spatial steps. We show the effect of the choice of $n_p$ on the
accuracy of the $y$-FT and the resulting convergence in $\Delta x \rightarrow 0$.}
\centering
\includegraphics[width=5.5in,height=2.5in]{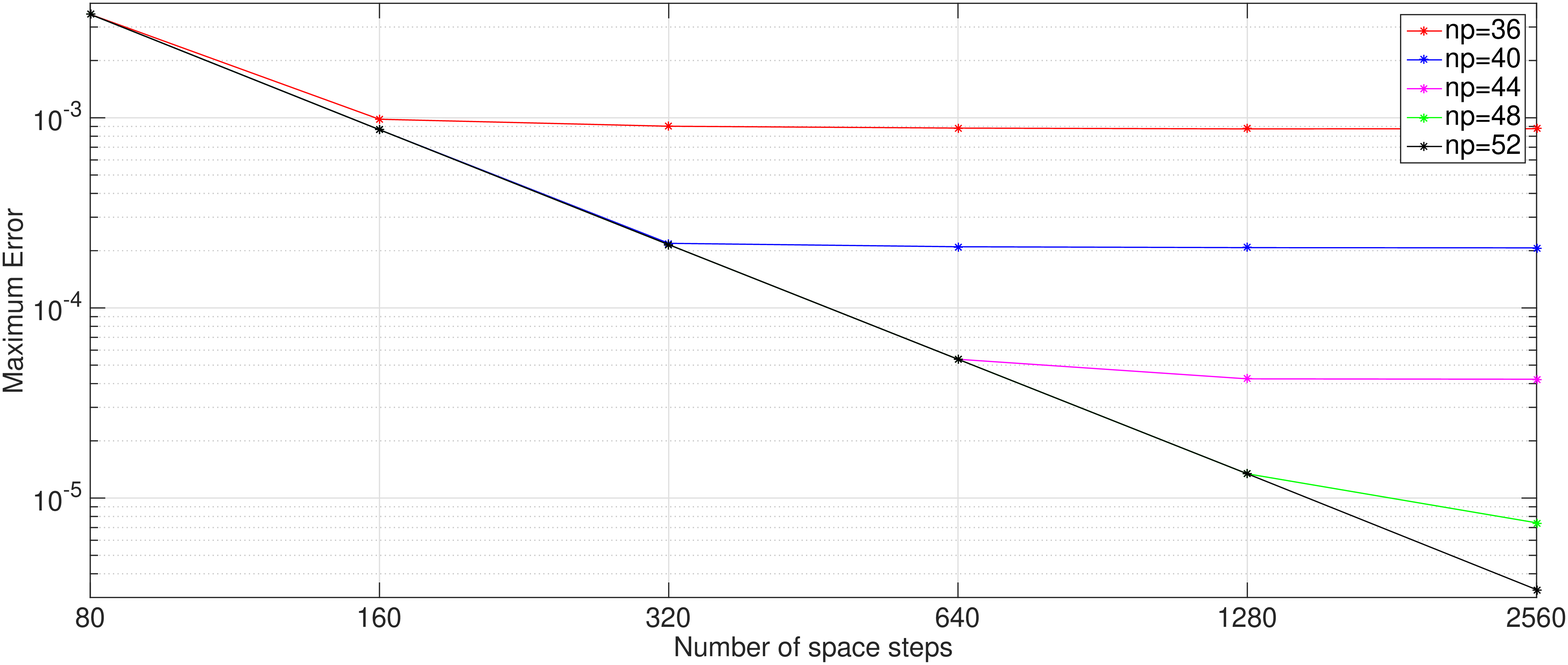}
\label{fig:SemiDLevelPlot}
\end{figure}

As $\Delta x$ is reduced we will find that we also need to reduce $\Delta p$ in order to maintain quadratic convergence. This behaviour is shown in
Figure \ref{fig:SemiDLevelPlot}.

Here, we have chosen the $x$-range as $[-10,10]$ and the $p$-range as $[-20,20]$ and we then monitor the convergence as
$\Delta x \rightarrow 0$ for a range of values of $\Delta p$, which correspond to $n_p=36, 40, 44, 48, 52$.

For large enough $n_p$, e.g.\ $n_p=52$, the plot shows quadratic convergence to the true solutions.
However, for a given value of $\Delta p$, quadratic convergence only continues to hold, as $\Delta x \rightarrow 0$, up to a critical value of $\Delta x$,
beyond which the decay of the discretization error stops.
Consequently, we need to reduce $\Delta p$ appropriately in order to achieve a satisfactory degree of accuracy.
We can also see that as $n_p$ increases towards the value where the error approaches the log-log line of quadratic convergence in $x$, the convergence in $p$ is very fast. Indeed, the horizontal asymptotes for large $n_x$, are roughly equally spaced on a log-scale when $n_p$ increases by a constant step, in line with the theoretical exponential convergence in $n_p$.

In the application in Section \ref{sec:hedge} we will study a setting where $u$ is interpreted as a joint probability density of two variables and we will be interested in the marginal density in $y$. We study this next.
The numerical solution, $U$, at time $T>0$, is calculated at grid points $(x_j,y_k)$ and is then numerically integrated over $x$ to give the marginal probability distribution for $y$.

\begin{figure}[H]
\caption{Toy model approximated by the semidiscrete Fourier method using the matrix exponential. Left, a comparison of marginal densities for $x\in[-20,20]$, $p\in[-20,20]$, $n_p=80$, $n_x=2560$. Right,  the difference between solutions.}
\centering
\includegraphics[width=3.4in,height=2in]{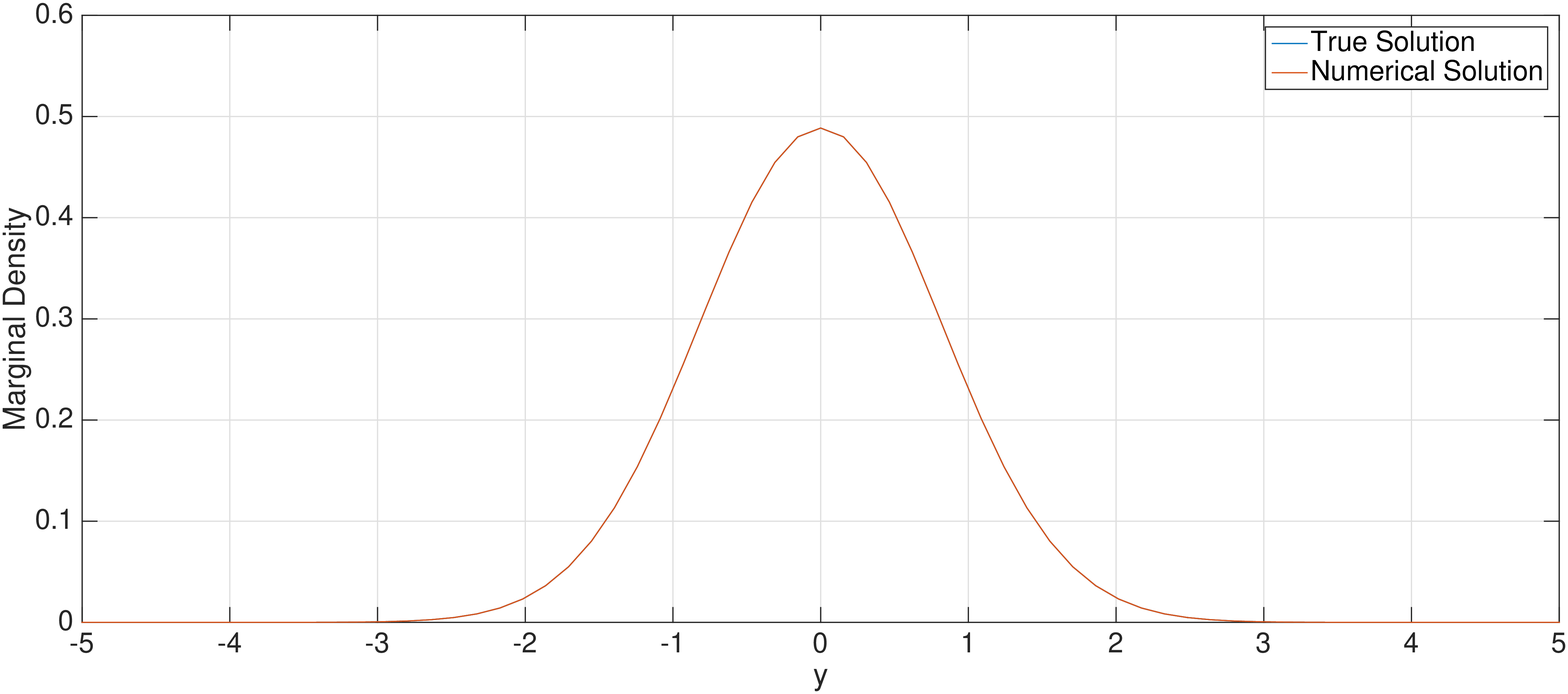} \hfill
\hspace{-0.5 cm}
\includegraphics[width=3.4in,height=2in]{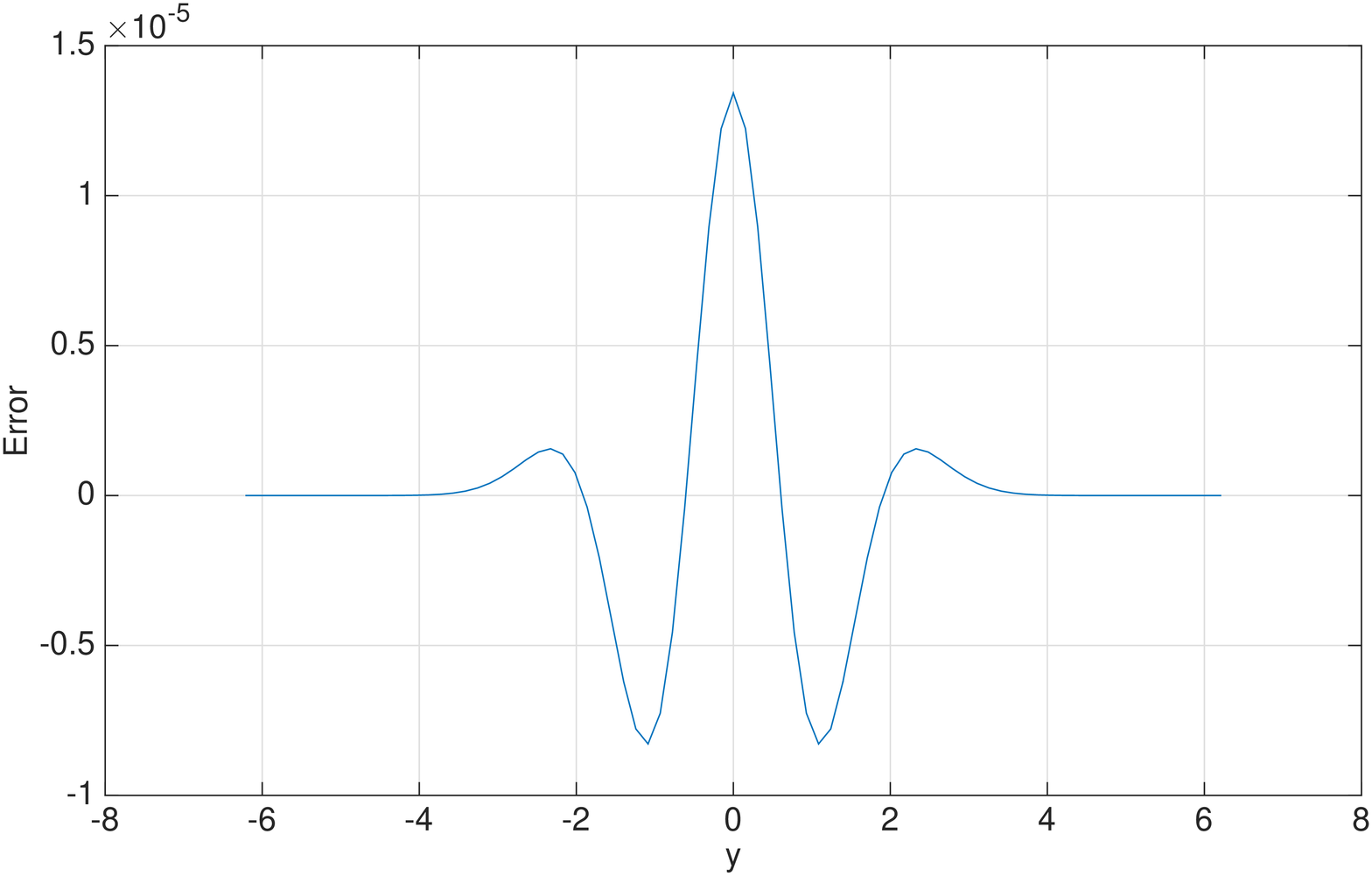}
\label{fig:SemiFComp2560}
\label{fig:SemiFDiff2560}
\end{figure}

Figure \ref{fig:SemiFComp2560}, left, plots the true and numerical marginal densities, while Figure \ref{fig:SemiFDiff2560}, right, shows the difference between the marginal
densities of the true and numerical solutions. We took $x\in[-20,20]$ and $p\in[-20,20]$ with 80 grid points in the $p$-direction and 2560 grid points in the $x$-direction.


%

%
%
%

\subsection{A fully-discrete Fourier scheme}
\label{subsec:fullydiscr}

The fully-discrete Fourier scheme again uses the $y$-FT to reduce the dimension of the problem, but it then applies discretization in both the
$x$- and $t$-direction.
After solving the linear system we again invert the $y$-FT to obtain
our approximate solution, $U(x_j,y,t_l)$, for the original problem. Here $t_l=l \Delta t$ are equally spaced time steps.
We use the toy model described above to study this scheme.


%
%

%
%

We discretize (\ref{toyftpde}) in the $x$-direction, using equally spaced
grid points, to obtain (\ref{eq:SEMIDSCHEME}), and then in the $t$-direction, using the (semi)implicit Euler scheme,
to obtain 
\begin{align*}
\frac{ V^{n+1}_{j}(p,t)- V^n_{j}(p,t)}{\Delta t}-ipx_jV^{n'}_{j}(p,t)=\frac{V^{n+1}_{j+1}(p,t)-2V^{n+1}_{j}(p,t)+V^{n+1}_{j-1}(p,t)}{\Delta x ^2} ,
\end{align*}
where $n'=n$ if the drift term is treated explicitly and $n'=n+1$ if it is treated implicitly. Obviously, other time discretizations are also possible, such
as second or higher order BDF, but for the sake of simplicity of the analysis to be presented we shall focus here on time stepping via the
Euler scheme.

\begin{remark}
For the purposes of an error analysis along the lines of Section \ref{sec:analysis}, we would apply the semidiscrete $x$-FT to obtain
\begin{equation*}
   \begin{aligned}
      \frac{W^{n+1}(s,p,t) -W^{n}(s,p,t) }{\Delta t} -p\frac{\partial W^{n'}(s,p,t) }{\partial s}  =-W^{n+1}(s,p,t) \frac{4\sin^2{(\frac{s\Delta x}{2})}}{\Delta x ^2},
   \end{aligned}
\end{equation*}
which can then be written as
\begin{equation}
   \begin{aligned}
      f(s)W^{n+1}(s,p,t)  -p\Delta t\frac{\partial W^{n'}(s,p,t) }{\partial s}  = W^{n}(s,p,t),
   \end{aligned}
\end{equation}
where $f(s)=1+4\lambda\sin^2{(\frac{s\Delta x}{2})}$ and $\lambda=\frac{\Delta t}{\Delta x^2}$. This equation is a `differential recursion' for the double-FT in the explicit case ($n'=n$) and an `integral recursion' in the implicit case ($n'=n+1$).
Although it is possible to identify leading order terms heuristically by an expansion, the rigorous analysis of these fully discrete schemes has proved elusive due to the complexity of the error propagation over the time steps across wave-numbers. We shall therefore explore these schemes by way of numerical tests.
\end{remark}




We present the results obtained by applying the two methods above to the toy model PDE. As with the semidiscrete Fourier method, the solution for $U$ is first integrated over $x$, before being plotted against $y$. We took $x\in[-20,20]$, $p\in[-20,20]$, $t \in [0,1]$, with $n_x=1800$ grid points in the $x$-direction, $n_p=80$ grid points in the $p$-direction, and $n_t=10240$ grid points in the $t$-direction.
We first present the results obtained by treating the drift term explicitly.



\begin{figure}[H]
\caption{Toy model approximated by the discrete Fourier method and an explicit approximation for the drift term. Left, the difference between the exact and approximate
solutions, with $x\in[-20,20]$, $p\in[-20,20]$, $t \in [0,1]$, $n_x=1800$, $n_p=80$, $n_t=10240$.
Right, with the other parameters the same, the error for decreasing timestep $\Delta t$ with $\Delta t/\Delta x^2$ held constant.}
\includegraphics[width=3.4in,height=2in]{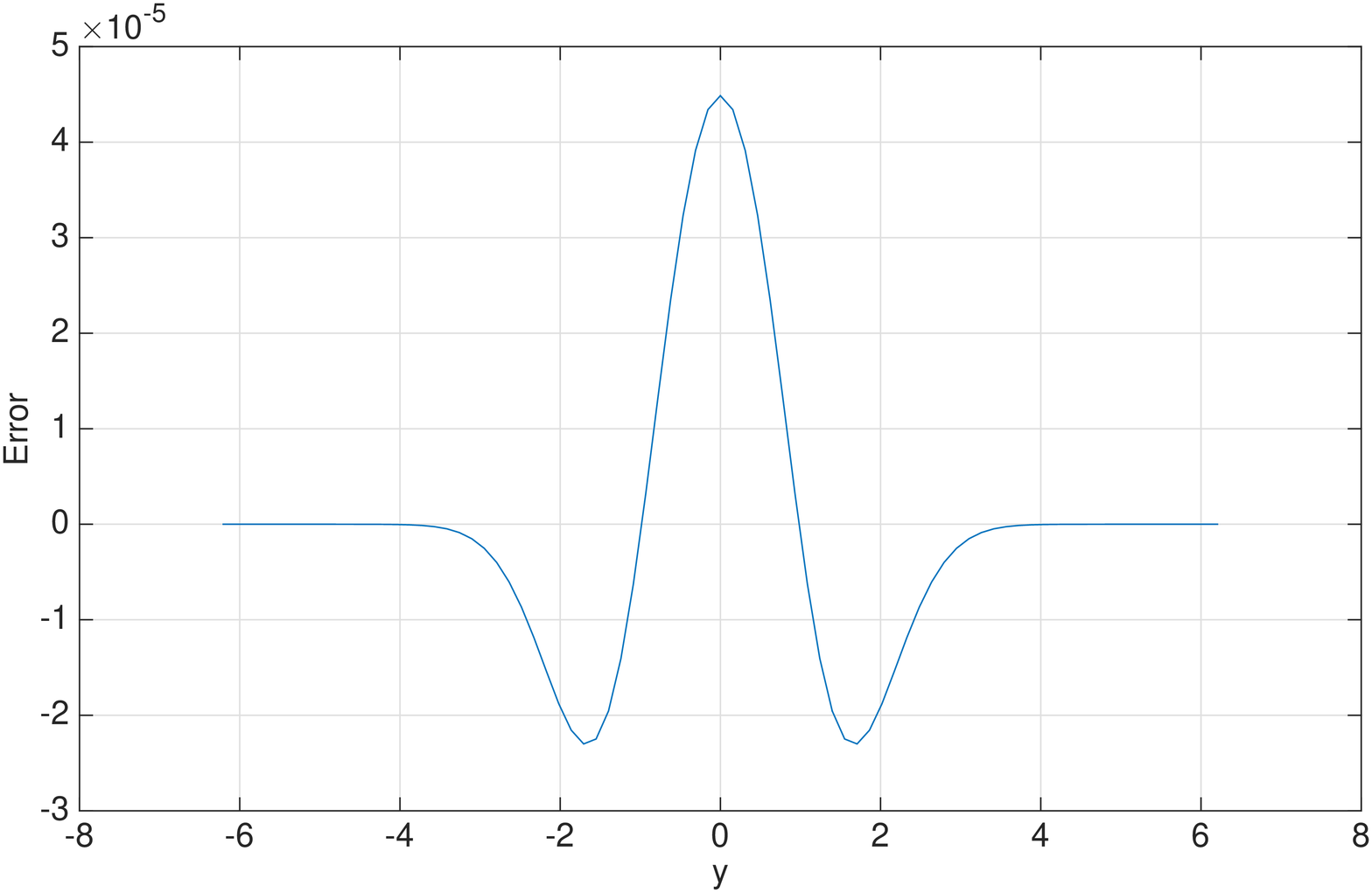} \hfill
\hspace{-0.5 cm}
\includegraphics[width=3.4in,height=2in]{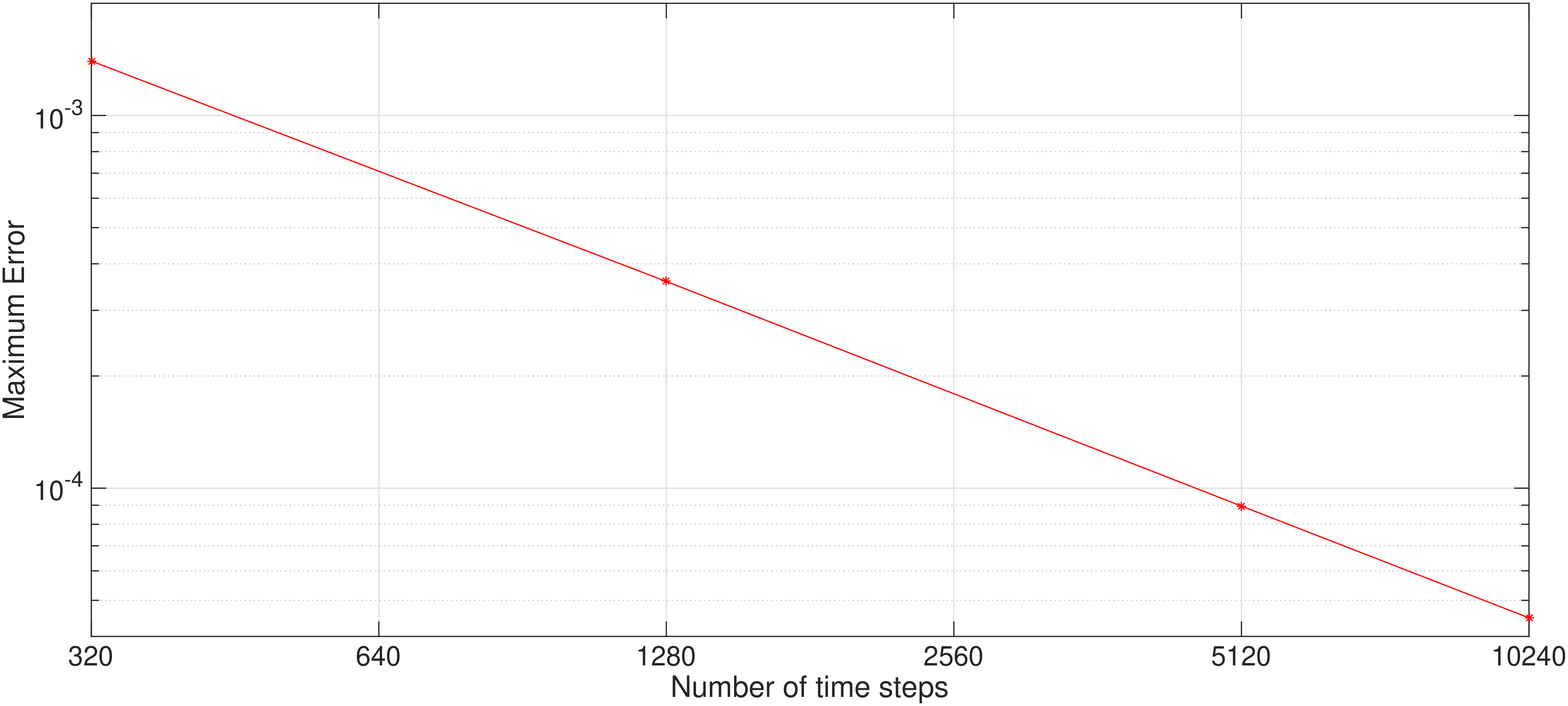}
\label{fig:DiffSolutionExplicitT10240}
\label{fig:ConvExplicitSep26}
\end{figure}

The left plot of Figure \ref{fig:DiffSolutionExplicitT10240} shows the difference between the true and numerical solutions. The numerical solution itself closely resembles the result in Figure \ref{fig:SemiFComp2560}, left. \\


%

Finally, the right plot of Figure \ref{fig:ConvExplicitSep26} shows the convergence rate for the solution as the grid is refined with $\lambda=\frac{\Delta t}{\Delta x^2}$ held constant.  We observe first-order convergence in time where the error is evaluated in the maximum norm. The slope of the last two points of the log-log plot is 1.00. \\


For this method, it is interesting to note that the maximum error in this case is approximately four times the corresponding error of the semidiscrete method (see Figures \ref{fig:SemiFDiff2560} and \ref{fig:DiffSolutionExplicitT10240}) and we observe upwards `bumps' at around $\pm 2$ in the latter case.

\begin{figure}[H]
\caption{Toy model approximated by the discrete Fourier method with an implicit approximation for the drift term. Left, the difference between the exact and approximate solutions,
with $x\in[-20,20]$, $p\in[-20,20]$, $t \in [0,1]$, $n_x=1800$, $n_p=80$, $n_t=10240$.
Right, with the other parameters the same, the error for decreasing timestep $\Delta t$ with $\Delta t/\Delta x^2$ held constant.
}
\includegraphics[width=3.4in,height=2in]{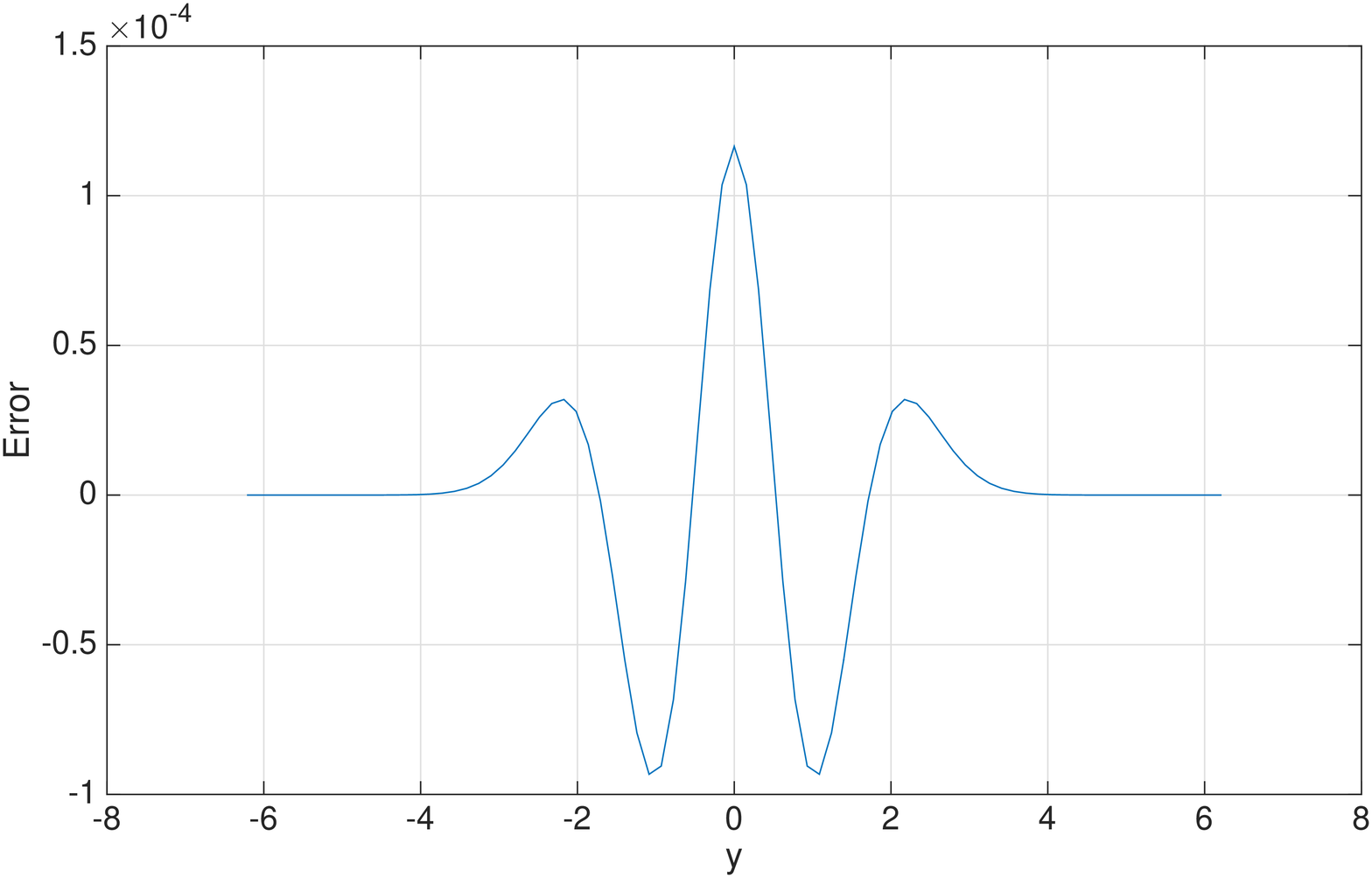}
\hspace{-0.5 cm}
\includegraphics[width=3.4in,height=2in]{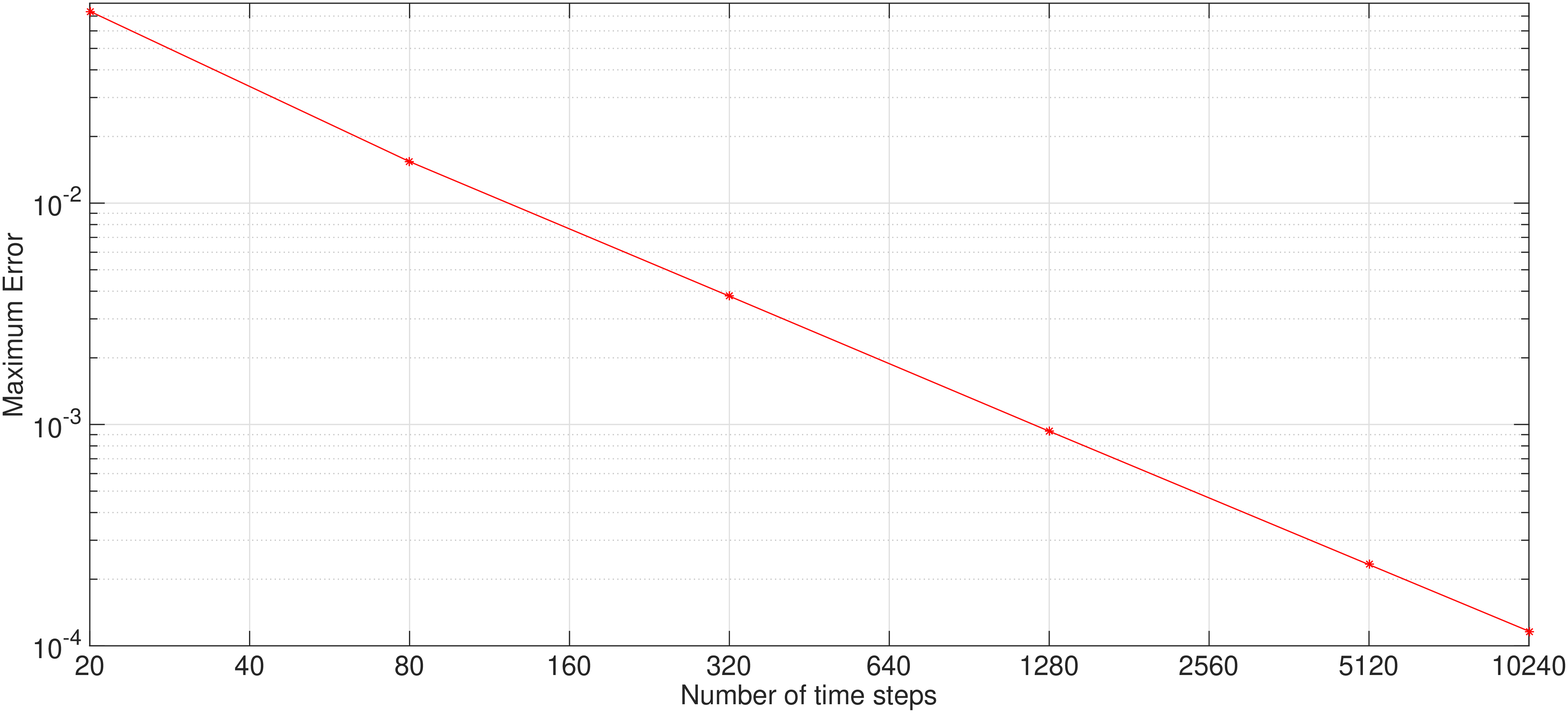}
\label{fig:ConvImplicitSep26}
\label{fig:DiffSolutionImplicitT10240}
\end{figure}

Next we present the results based on treating the drift implicitly.
%
%
The left plot in Figure \ref{fig:DiffSolutionImplicitT10240} shows the difference between the true and numerical solutions. Once again, the numerical solution itself closely resembles the result in Figure \ref{fig:SemiFComp2560}. \\


%

Finally, the right plot in Figure \ref{fig:ConvImplicitSep26} shows the convergence rate for the solution as the grid is refined. We observe first-order convergence in time where the error is evaluated in the maximum norm. The slope of the last two points of the log-log plot is 1.00. \\


We observe that the maximum error in this case is approximately twice the corresponding error of the explicit drift case (see Figures \ref{fig:DiffSolutionExplicitT10240} and \ref{fig:DiffSolutionImplicitT10240}).

\section{Application to a financial hedging problem}
\label{sec:hedge}

In this section, we apply the Fourier method to a hypoelliptic problem from mathematical finance, namely the computation of hedging errors under misspecification of the market model. If the true model governing an underlying stock is known to the trader and the market is ``complete'', they can perfectly hedge a position in an option by dynamic trading in the stock. If the model is unknown and the hedging strategy is based on a misspecified model, a certain profit or loss will materialize and one can ask what the distribution of this hedging error is.

Our example follows a simplified version of the analysis in El Karoui \emph{et al}., \cite{ELK}, where the true dynamics of a single underlying asset satisfy the SDE
\begin{equation}
\frac{\dd S_t}{S_t}=\mu_t \dd t + \sigma_t \dd W_t,
\end{equation}
where $\sigma$ is a volatility process, $\mu$ a drift, and $W$ a standard Brownian motion. 
We analyze the situation where the trader instead assumes
that the underlying follows a different process, i.e. mis-specifies the model for $S$ as
\begin{equation} \label{eq:DiffVol}
\frac{\dd S_t}{S_t}=\widehat{\mu}_t \dd t + \widehat{\sigma} \dd W_t
\end{equation}
for some 
constant $\widehat{\sigma}$.
We define $\widehat{V}(S_t,t)$ to be the (Black-Scholes) price of the European option with payoff $F(S_T)$ based on the process defined in (\ref{eq:DiffVol}), which satisfies the PDE
\begin{eqnarray}
\label{bspde}
\mathcal{L}^{\widehat{\sigma}}\widehat{V} \equiv \frac{\partial \widehat{V}}{\partial t}+\frac{1}{2}\widehat{\sigma}^2S^2 \frac{\partial^2 \widehat{V}}{\partial S^2}+r S\frac{\partial \widehat{V}}{\partial S}-r\widehat{V} = 0.
\end{eqnarray}

If the trader then uses the sensitivity $\frac{\partial \widehat{V}}{\partial S}$ of $\widehat{V}$ as the hedge ratio,
it is shown in \cite{ELK} that the hedging ``error'' $Y_t$, i.e, the difference in time $t$ value between the option and the portfolio set up to hedge it,
discounted with the risk-free rate $r$ to the present time,
is governed by the stochastic differential equation
%
\begin{equation}
\dd Y_t= \frac{1}{2} \mbox{e}^{-rt}(\widehat{\sigma}^2-\sigma_t^2)S_t^2 \frac{\partial^2\widehat{V}}{\partial S^2}\dd t.
\end{equation}
%
%
%
We notice that the dynamics of $Y_t$ only consists of a drift term, i.e.\ there is no Brownian component.
In the following, we assume that $\sigma_t=\sigma$ constant, but different from $\widehat{\sigma}$.
%
%

Therefore, if we write down the Kolmogorov Forward Equation (KFE) for the probability density function $P$ of the pair $(S_t,Y_t)$ at time $t$,
we find
\begin{eqnarray}\nonumber
\frac{\partial P}{\partial t} &=&
\frac{1}{2} \sigma^2 \frac{\partial^2}{\partial S^2} \left(S^2 P\right)
- \mu  \frac{\partial}{\partial S} \left(S P\right) - c \frac{\partial P}{\partial y} \\
&=& \frac{1}{2}\sigma^2 S^2 \frac{\partial^2P}{\partial S^2}+(2\sigma^2-\mu) S \frac{\partial P}{\partial S}+(\sigma^2-\mu)P -c\frac{\partial P}{\partial y},
\label{eq:KFEBlackScholes}
\end{eqnarray}
where
\begin{eqnarray}\label{eq:KFEBlackScholesDrift}
c(S,t)=\frac{1}{2}\mbox{e}^{-rt}(\widehat{\sigma}^2-\sigma^2)S^2 \frac{\partial^2\widehat{V}}{\partial S^2}.
\end{eqnarray}
The equation hence falls into the class of hypoelliptic PDEs studied in the preceding sections.
%



%
%
%

We solve the KFE numerically using the following parameters, $\sigma=0.1$, $\widehat{\sigma}=0.2$, $\mu=0.05$, and $r=0.05$. We take $T=1.0$, $S_0=125$, and the option being hedged is a so-called spread with payoff defined as
\begin{align*}
F(S_T)=\max(S_T-100,0)-\max(S_T-150,0).
\end{align*}
We note that the convexity of this option payoff changes as $S$ varies, and hence also the sign of the function $c$, which contains the second derivative of $\widehat{V}$. The function $\widehat{V}$ is known in closed form in this case as the solution to (\ref{bspde}), and so is its second derivative, referred to as the ``gamma'' in the financial community.

To find the univariate density of the discounted terminal hedging error, the joint density is integrated over the $S$-direction.
The univariate distribution of $y$ has been calculated using the fully-discrete Fourier method from Section \ref{subsec:fullydiscr}.
A Monte Carlo hedging simulation has been also performed where the asset price was simulated according to the true asset model (which used the true volatility) while the hedged portfolio was controlled by delta-hedging using the Black--Scholes equation with the erroneous volatility.

\begin{figure}[H]
\caption{Black--Scholes model with mis-specified volatility. Comparison of the KFE solution, integrated over $x$, for an increasing number of spatial steps $n_S$ and keeping $n_t/n_S^2$ fixed. Shown also is the empirical density from a Monte Carlo simulation as detailed in the text.
}
\begin{center}
\includegraphics[width=5.5in,height=2.5in]{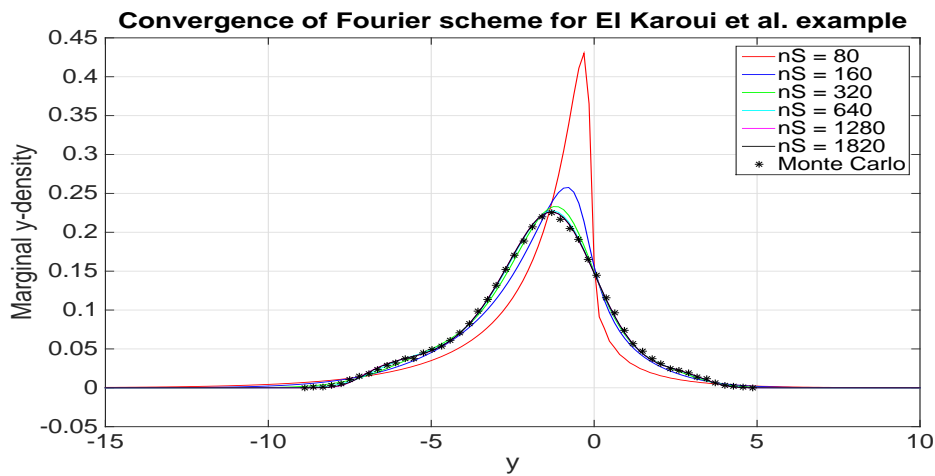}
\end{center}
\label{fig:ElKarouiCompNewNov2015}
\end{figure}

Figure~\ref{fig:ElKarouiCompNewNov2015} shows how the approximate distributions compare.
The shape of the hedging error distribution can be rationalized by noting that the positive drift, $\mu$, implies that the stock price is more likely to move into a regime where the option gamma is negative, resulting in a preponderance of negative hedging errors. Note that in these calculations, we assume that the option is sold for a price consistent with $\widehat{\sigma}$ so that $Y$ starts from zero. If some other price were initially realized, it would result in the densities being shifted by this amount.
%
%

For the fully-discrete Fourier method, we used $S \in [5,245]$, $y \in [-20,20]$ (not plotted over whole range) and $p \in [-20,20]$, $n_p=n_y=220$. For the Monte Carlo results, we used 100 time steps and 100000 paths as well as 50 bins for the approximation of the density.

Having derived an approximation to the hedging error distribution, one can compute quantities of financial interest, such as the Value-at-Risk, a quantile-based risk-measure. For instance, the 10\%-quantile of the hedging error for the data set above is found to be $-4.58$. This compares to an initial option value of 29.31 and 26.67, respectively, for the low and high volatility.

\section{Multi-dimensional problems}
\label{ch:OTHER}
\label{sec:multi}

It is straightforward to extend the toy model (\ref{toypde}), (\ref{toydelta}) to higher dimensions by introducing multiple diffusive terms and/or multiple drift terms.
For the semidiscrete case, we expect that the analysis will be similar to that carried out for the one-dimensional case. We present this analysis in the remainder of this section.

\subsection{Multiple diffusive terms} \label{app:HIGHERDIMTOY}


We consider the case of a PDE with a single drift coefficient which is a linear combination of the components of $x$, 
\begin{alignat}{2}
\label{general-2d}
\frac{\partial u}{\partial t}+(\gamma_1 x_1 +\gamma_2 x_2 ) \frac{\partial u}{\partial y} &=
a_{11} \frac{\partial^2 u}{\partial x_1^2} + 2 a_{12} \frac{\partial^2 u}{\partial x_1 \partial x_2} +a_{22} \frac{\partial^2 u}{\partial x_2^2},
&&\quad (x_1,x_2,y,t) \in \mathbb{R}^3
\times (0,T], \\
u(x_1,x_2,y,0)&=\delta (x_1-x_{1,0})\otimes \delta (x_2-x_{2,0})\otimes \delta (y-y_0), &&\quad (x_1,x_2,y) \in \mathbb{R}^3, 
\nonumber
\end{alignat}
which extends the toy model to one containing two independent variables $x_1$ and $x_2$.
We consider the hypoelliptic case, i.e., where $(\gamma_1,\gamma_2) \neq (0,0)$ and $(a_{ij})_{1\le i,j,\le 2}$ is strictly positive definite,
i.e., $a_{11}, a_{22}>0$ and $a_{11} a_{22}> a_{12}^2$.

We apply the $y$-FT to get
\begin{eqnarray}
\label{2dft}
\frac{\partial v}{\partial t}-ip(\gamma_1 x_1 +\gamma_2 x_2) v =a_{11} \frac{\partial^2 v}{\partial x_1^2}
+ 2a_{12} \frac{\partial^2 v}{\partial x_1 \partial x_2}+a_{22} \frac{\partial^2 v}{\partial x_2^2}, \quad (x_1,x_2,p,t) \in \mathbb{R}^3 \times (0,T].
\end{eqnarray}
Now for the analytical solution we next apply the $x_1$-FT and the $x_2$-FT in turn to get
\begin{eqnarray}
\label{true-double-ft}
\frac{\partial w}{\partial t}-p \bigg(\gamma_1 \frac{\partial w}{\partial s_1} +\gamma_2  \frac{\partial w}{\partial s_2} \bigg) =-(a_{11} s_1^2+2a_{12} s_1s_2+a_{22} s_2^2)w,\quad (s_1,s_2,p,t) \in \mathbb{R}^3 \times (0,T],
\end{eqnarray}
in conjunction with the initial condition
\[
w(s_1,s_2,p,0)=1.
\]
We can now solve this first-order hyperbolic initial-value problem directly, but we can also use an ansatz by seeking a solution of the form
\[
w(s_1,s_2,p,t)=\exp{(-a_{11}s_1^2t -2a_{12}s_1s_2t -a_{22}s_2^2t xs- B_1s_1pt^2-B_2s_2pt^2-Cp^2t^3)}
\]
and then finding the coefficients that fit the PDE.
%
%
%
%
%
%
%
%
By insertion we get
\begin{align*}
B_1 &= \gamma_1a_{11}+\gamma_2a_{12},  \\
B_2 &= \gamma_1a_{12}+\gamma_2a_{22},  \\
C &= \frac{1}{3}( \gamma_1(a_{11} +a_{12}  )+\gamma_2(a_{12} +a_{22}) ),
\end{align*}
and thus we have the analytical solution in Fourier space. We could, if necessary, invert all these transforms but we will instead concentrate on the effect of discretization in the $x_1$- and $x_2$-directions.
So we now start with the PDE (\ref{2dft})
%
%
and discretize this as follows by writing
%
\begin{eqnarray*}
\frac{\partial V_{j,k}}{\partial t}-ip(\gamma_1 x_{1,j} +\gamma_2 x_{2,k})V_{j,k} &=&
 a_{11} \frac{V_{j+1,k}-2V_{j,k}+V_{j-1,k}}{\Delta x_1^2} +  \, a_{22}  \frac{V_{j,k+1}-2V_{j,k}+V_{j,k-1}}{\Delta x_2^2} \\
&&\qquad+\, 2a_{12}  \frac{V_{j+1,k+1}-V_{j+1,k-1}-V_{j-1,k+1}+V_{j-1,k-1} }{4\Delta x_1\Delta x_2},
\end{eqnarray*}
where we have omitted $p$ and $t$ as arguments for the sake of brevity, and where $V_{j,k}(p,t)$ is an approximation to $v(x_{1,j},x_{2,k},p,t)$ on a uniform
two-dimensional mesh of widths $\Delta x_1, \Delta x_2>0$.
We can then apply the Fourier transforms in the $x_1$ and $x_2$ variables and after steps similar to the one-dimensional case
obtain the solution for the double-transformed equation,
\begin{equation}
W(s_1,s_2,p,t)
=
w(s_1,s_2,p,t)
\exp{\bigg(\frac{1}{p}  \int^{s_1+pt}_{s_1} ( a_{11}g_{11} (\sigma,s_2)+2a_{12}g_{12} (\sigma,s_2)+a_{22}g_{22} (\sigma,s_2)     )\dd\sigma \bigg)},
\end{equation}
where
\begin{align*}
g_{11} (s_1,s_2) &= s_1^2-\frac{4}{\Delta x_1^2}\sin^2{\bigg(\frac{s_1\Delta x_1}{2}\bigg)},  \\
g_{12} (s_1,s_2) &= s_1s_2- \frac{1}{\Delta x_1\Delta x_2}\sin{(s_1\Delta x_1)}\sin{(s_2\Delta x_2)}, \\
g_{22} (s_1,s_2) &= s_2^2-\frac{4}{\Delta x_2^2}\sin^2{\bigg(\frac{s_2\Delta x_2}{2}\bigg)}.
\end{align*}
We can now investigate the low and high wavenumber behaviour in terms of the new variables.
%
%
For the sake of simplicity of the exposition, we shall do this for the special case $\gamma_1 = 1$, $\gamma_2 = 0$, $a_{11} = a_{22} = 1$, $a_{12} = \rho \in (-1,1)$, i.e., we consider
\begin{eqnarray}
\label{special-2d}
\frac{\partial u}{\partial t} + x_1 \frac{\partial u}{\partial y} =  \frac{\partial^2 u}{\partial x_1^2} + 2 \rho \frac{\partial^2 u}{\partial x_1 \partial x_2} +
\frac{\partial^2 u}{\partial x_2^2}.
\end{eqnarray}

\begin{remark}
Equation (\ref{special-2d}) is a simplified case of the general equation (\ref{general-2d}), where the term $x_2$ has been dropped
from the drift and the diffusion is normalized. This can be achieved by a rotation of the original co-ordinates to align the drift
$(\gamma_1,\gamma_2)$ with the $x_1$-axis, and by subsequent scaling of the new $y$, $x_1$ and $x_2$ co-ordinates.
\end{remark}

In this case, the numerical solution is
\begin{eqnarray*}
W(s_1,s_2,p,t)
=w(s_1,s_2,p,t)
\exp{\bigg(\frac{1}{p}  \int^{s_1+p t}_{s_1} (g_{11} (\sigma,s_2)+ 2 \rho g_{12} (\sigma,s_2) + g_{22} (\sigma,s_2) )
\dd\sigma \bigg)}.
\end{eqnarray*}
Similarly to the one-dimensional case, this leads us to investigate
\[
W(s_1,s_2,p,t)
=
W_1(s_1,p,t) \ \exp \bigg(- \frac{1}{p}  \int^{s_1+p t}_{s_1}
\left(
2 \rho \frac{\sin{(\sigma\Delta x_1)}\sin{(s_2\Delta x_2)}}{\Delta x_1\Delta x_2}
+ \frac{4}{\Delta x_2^2}\sin^2{\bigg(\frac{s_2\Delta x_2}{2}\bigg)} \right) \dd\sigma \bigg),
\]
where $W_1$ is the solution from the one-dimensional case.
By standard integration and the change of variables $\xi = pt/2$ and $\eta_1 = s_1 + pt/2$ we find this to be
\begin{eqnarray*}
W(\eta_1,s_2, \xi,t) = W_1(\eta_1,\xi,t) \, \exp \bigg(- 2 t \rho \sinc(\xi \Delta x_1) \frac{\sin(\eta_1 \Delta x_1)}{\Delta x_1}  \frac{\sin(s_2 \Delta x_2)}{\Delta x_2}  - \frac{4t}{\Delta x_2^2} \sin^2(s_2 \Delta x_2/2) \bigg),
\end{eqnarray*}
where
\begin{eqnarray*}
W_1(\eta_1,\xi,t) =
\exp \bigg(-
\frac{2t}{\Delta x_1^2} \left(1-\sinc(\xi \Delta x_1) \cos(\eta_1 \Delta x_1)\right) \bigg).
\end{eqnarray*}
The exact solution in these variables is
\begin{eqnarray}
\label{2d-exact}
w(\eta_1,s_2,\xi,t)=\exp{\left(- t \left\{ \eta_1^2 + \frac{\xi^2}{3} + s_2^2 + 2 \rho \eta_1 s_2 \right\}\right)}
\end{eqnarray}
and we recognize in the first two terms the one-dimensional solution
\begin{eqnarray}
\label{1d-exact}
w_1(\eta_1,\xi,t)=\exp{\left(- t \left\{ \eta_1^2 + \frac{\xi^2}{3} \right\}\right)}.
\end{eqnarray}
We proceed by a wavenumber analysis broadly similar to the one before, but made somewhat more complicated by the
presence of an extra variable $s_2$ and the fact that the problem degenerates as $|\rho|\rightarrow 1$, necessitating a
more careful estimation for $\rho$ close to 1.

\subsection{Joint low wavenumbers}

We write
\begin{eqnarray*}
\log W = \log W_1
- t  \left\{
2  \rho  \sinc(\xi \Delta x_1) \frac{\sin(\eta_1 \Delta x_1)}{\Delta x_1} \frac{\sin\left(s_2 \Delta x_2\right)}{ \Delta x_2}
+ \frac{4}{\Delta x_2^2} \sin^2\left(\frac{s_2 \Delta x_2}{2}\right)
\right\}
\end{eqnarray*}
and compare this to
\begin{eqnarray*}
\log w = \log w_1 - t \{s_2^2 + 2 \rho \eta_1 s_2 \},
\end{eqnarray*}
where $w$ is the exact solution to the two-dimensional problem from (\ref{2d-exact}) and $w_1$ is the solution to the one-dimensional problem
from (\ref{1d-exact}). By standard Taylor expansion we find
\begin{align*}
\log(W/w) & = \log(W_1/w_1) + 2 \rho t \left(
\frac{1}{12} \eta_1 s_2^3 \Delta x_2^2 + \frac{1}{6} \eta_1^3 s_2^2 \Delta x_1^2 + \eta_1 s_2^2 \xi^2 \Delta x_1^2
\right) + \frac{t}{12} s_2^4 \Delta x_2^2 + o(\Delta x_1^2) + o(\Delta x_2^2) \\
&= \log(W_1/w_1) + \left(\frac{\rho t}{3} \eta_1^3 s_2^2 + 2 \rho t \eta_1 s_2^2 \xi^2 \right) \Delta x_1^2 + \left(\frac{\rho t}{6} \eta_1 s_2^3 + \frac{t}{12} s_2^4 \right) \Delta x_2^2 + o(\Delta x_1^2) + o(\Delta x_2^2)\\
&=
 \left(
 \frac{ 2t}{4!} \eta_1^4  + \frac{2 t}{5!} \xi^4 +
 \frac{\rho t}{3} \eta_1^3 s_2^2 + 2 \rho t \eta_1 s_2^2 \xi^2 \right) \Delta x_1^2 + \left(\frac{\rho t}{6} \eta_1 s_2^3 + \frac{t}{12} s_2^4 \right) \Delta x_2^2
 + o(\Delta x_1^2) + o(\Delta x_2^2).
\end{align*}
Therefore the numerical error contribution from the low wavenumber regime is
\begin{align}
\label{2derr}
I_1(x_1,x_2,y,t) &=
\Delta x_1^2 F_1(x_1,x_2,y,t)
+
\Delta x_2^2 F_2(x_1,x_2,y,t)  + o(\Delta x_1^2) + o(\Delta x_2^2),
\end{align}
where
\begin{align*}
F_1 &=
t^2 \left(
 \frac{1}{4!} \left(\frac{\partial }{\partial x_1} + \frac{t}{2} \frac{\partial }{\partial y} \right)^4  + \frac{1}{5!}  \left(\frac{t}{2} \frac{\partial }{\partial y} \right)^4  \right) \ u\\
& \hspace{2 cm} + \ t^2 \left(
 \frac{\rho}{3!}  \left(\frac{\partial }{\partial x_1} + \frac{t}{2} \frac{\partial }{\partial y} \right)^3 \frac{\partial^2 }{\partial x_2^2} +
 \rho \left(\frac{\partial }{\partial x_1} + \frac{t}{2} \frac{\partial }{\partial y} \right) \frac{\partial^2 }{\partial x_2^2} \left(\frac{t}{2} \frac{\partial }{\partial y} \right)^2
 \right) \ u \\
F_2 &= t^2 \left(
\frac{2 \rho}{4!}  \left(\frac{\partial }{\partial x_1} + \frac{t}{2} \frac{\partial }{\partial y} \right) \frac{\partial^3 }{\partial x_2^3}
+ \frac{1}{4!} \frac{\partial^4 }{\partial x_2^4}
\right) \ u.
\end{align*}
The point is less the explicit form, but that the error can be expressed in terms of up to fifth mixed partial derivatives, or
up to fourth if $\rho=0$. As in Section \ref{sec:analysis}, we will find that this is the only wavenumber range that contributes to the leading order error, and therefore, up to higher order terms, (\ref{2derr}) fully describes the discretization error.

\subsection{Low $\boldsymbol{\xi}$-wavenumbers and high $\boldsymbol{\eta_1}$- or $\bf s_2$-wavenumbers}

We require two simple inequalities. The first one states that for $|\eta_1 \Delta x_1| \le \pi/2$, and since
$\sinc(\xi \Delta x_1)>0$ (in the present small $\xi$ regime),
\[
1 - \sinc(\xi \Delta x_1) \cos(\eta_1 \Delta x_1) \ge 1 - \cos(\eta_1 \Delta x_1) = 2 \sin^2(\eta_1 \Delta x_1/2).\]
The second elementary inequality used in the argument below, in the transition from the second-to-last inequality to the last
inequality, is that
\[
\sin \alpha \ge \frac{2}{\pi} \alpha\qquad \mbox{for}\qquad 0\le \alpha \le \frac{\pi}{2}.
\]
Therefore
\begin{align*}
W(\eta_1,s_2) &\le\exp \bigg(
- \frac{4 t}{\Delta x_1^2} \sin^2(\eta_1 \Delta x_1/2) + 2 t |\rho|  \frac{|\sin(\eta_1 \Delta x_1)|}{\Delta x_1}
\frac{|\sin(s_2 \Delta x_2)|}{\Delta x_2} - \frac{4 t}{\Delta x_2^2} \sin^2(s_2 \Delta x_2/2)
\bigg) \\
 &\le
\exp \bigg(
- \frac{4 t}{\Delta x_1^2} \sin^2(\eta_1 \Delta x_1/2) + 8 t |\rho|  \frac{|\sin(\eta_1 \Delta x_1/2)|}{\Delta x_1}
\frac{|\sin(s_2 \Delta x_2/2)|}{\Delta x_2} - \frac{4 t}{\Delta x_2^2} \sin^2(s_2 \Delta x_2/2)
\bigg) \\
&=
\exp \bigg(
- 4 t (1-|\rho|)
\left(
\frac{1}{\Delta x_1^2} \sin^2(\eta_1 \Delta x_1/2) + \frac{1}{\Delta x_2^2} \sin^2(s_2 \Delta x_2/2)
\right)
\bigg)\\
&
\hspace{5 cm} \cdot \; \exp \bigg(
- 4t |\rho| \left(
\frac{\sin(\eta_1 \Delta x_1/2) }{\Delta x_1} - \frac{\sin^2(s_2 \Delta x_2)}{\Delta x_2}
\right)^2
\bigg) \\
&\leq \exp \bigg(- 4 t (1-|\rho|) \left(
\frac{1}{\Delta x_1^2} \sin^2(\eta_1 \Delta x_1/2) + \frac{1}{\Delta x_2^2}
\sin^2(s_2 \Delta x_2/2)
\right)
\bigg) \\
&\le
\exp \bigg(
-\frac{4 \, t \, (1-|\rho|)}{\pi^2}
\left(\eta_1^2 + s_2^2 \right) \bigg) \\
&= o(\Delta x^r)\qquad \mbox{as $\Delta x \rightarrow 0$}
\end{align*}
for any $r>0$ if either $\nu_1 \ge \Delta x_1^{-q}$ or $s_2\ge \Delta x_2^{-q}$ for any $q>0$.

For $|\eta_1 \Delta x_1| > \pi/2$, we have that
\[
1 - \sinc(\xi \Delta x_1) \cos(\eta_1 \Delta x_1) \ge 1,
\]
and therefore
\begin{align*}
W(\eta_1,s_2)
&\le
\exp \bigg(
- \frac{2 t}{\Delta x_1^2} + 2 t |\rho|  \frac{1}{\Delta x_1}
\frac{|\sin(s_2 \Delta x_2)|}{\Delta x_2} - \frac{4 t}{\Delta x_2^2} \sin^2(s_2 \Delta x_2/2)
\bigg) \\
&\le
\exp \bigg(
- \frac{2 t}{\Delta x_1^2}  + \frac{4 t}{\Delta x_1}
\frac{|\sin(s_2 \Delta x_2/2)|}{\Delta x_2} - \frac{4 t}{\Delta x_2^2} \sin^2(s_2 \Delta x_2/2)
\bigg) \\
&\le
\exp \bigg(
- \frac{2 t}{\Delta x_1^2}  + \frac{\sqrt{21} t}{\Delta x_1}
\frac{|\sin(s_2 \Delta x_2/2)|}{\Delta x_2} - \frac{4 t}{\Delta x_2^2} \sin^2(s_2 \Delta x_2/2)
\bigg) \\
&=
\exp \bigg(
-  \frac{t}{2}
\left(
\frac{1}{\Delta x_1^2} + \frac{1}{\Delta x_2^2} \sin^2(s_2 \Delta x_2/2)
\right)
\bigg) \cdot \\
&
\hspace{5 cm} \cdot \; \exp \bigg(
- t \left(
\sqrt{\frac{3}{2}} \frac{1}{\Delta x_1} - \sqrt{\frac{7}{2}} \frac{|\sin(s_2 \Delta x_2/2)|}{\Delta x_2}
\right)^2
\bigg) \\
&\le
\exp \bigg(
-  \frac{t}{2}
\left(
\frac{1}{\Delta x_1^2} + \frac{1}{\Delta x_2^2} \sin^2(s_2 \Delta x_2/2)
\right)
\bigg) \\
& = o(\Delta x^r)
\end{align*}
as before.

\subsection{High $\boldsymbol{\xi}$-wavenumbers}

Completing squares, we write
\begin{align*}
W(\eta_1,s_2,\xi) &=
W_1(\eta_1,\xi) \, \cdot \,
\exp\left(
\frac{t}{\Delta x_1^2} \cos^2(s_2\Delta x_2/2) \sinc^2(\xi \Delta x_1) \sin^2(\eta_1 \Delta x_1) \rho^2
\right) \, \cdot \\
&\qquad \cdot \, \exp\left(
- t \left(
\frac{\rho}{\Delta x_1}
\cos(s_2 \Delta x_2/2) \sinc(\xi \Delta x_1) \sin(\eta_1 \Delta x_1)
+ \frac{2}{\Delta x_2}
\sin(s_2 \Delta x_2/2)
\right)^2
\right),
\end{align*}
and so, neglecting the last factor, which is $\le 1$,
\begin{align*}
\log W &\le - \frac{2 t}{\Delta x_1^2}
\left(1-\sinc(\xi \Delta x_1) \cos(\eta_1 \Delta x_1)\right)
+ \rho^2 \frac{t}{\Delta x_1^2} \cos^2(s_2\Delta x_2/2) \sinc^2(\xi \Delta x_1) \sin^2(\eta_1 \Delta x_1) \\
&\le \frac{t}{\Delta x_1^2}
\left(-2 + 2 \sinc(\xi \Delta x_1) \cos(\eta_1 \Delta x_1) + \sinc^2(\xi \Delta x_1) \sin^2(\eta_1 \Delta x_1) \right) \\
&\le \frac{t}{\Delta x_1^2} \left(-2 + |\sinc(\xi \Delta x_1)|
\left(2 |\cos(\eta_1 \Delta x_1)| + \sin^2(\eta_1 \Delta x_1) \right)  \right) \\
&= \frac{t}{\Delta x_1^2} \left(-2 + |\sinc(\xi \Delta x_1)|
\left(2 - (1-|\cos(\eta_1 \Delta x_1)|)^2 \right)
\right) \\
&\le - \frac{2 t}{\Delta x_1^2} (1-|\sinc(\xi \Delta x_1)|),
\end{align*}
and following the remainder of the argument in the one-dimensional case we can deduce that the contribution from
this range is also exponentially small.

%
%
%
%
%

\subsection{Multiple drift terms} \label{app:MOREDRIFTTOY}

Another possible extension is a model of the form
\begin{equation} \label{eq:MULTDRIFT}
\frac{\partial u}{\partial t} +a_1 x \frac{\partial u}{\partial y_1}  +a_2 x \frac{\partial u}{\partial y_2} = \frac{\partial^2 u}{\partial x^2},
\end{equation}
where $(a_1,a_2) \neq  (0,0)$ so that we have two drift terms and neither of the drift coefficients depend on $y$.
%
%
%
%
%
%
%
%
%
For this, we can apply the $y_{1}$-FT and the $y_{2}$-FT to get
\[ \frac{\partial v}{\partial t}-ia_1 p_1xv-ia_2 p_2xv= \frac{\partial^2 v}{\partial x^2},   \]
where $p_1$ and $p_2$ are wavenumbers.
To find the analytical solution we apply the $x$-FT to get
%
%
%
\[ \frac{\partial w}{\partial t} -(a_1p_1+a_2 p_2)\frac{\partial w}{\partial s}= -s^2w  .\]
The analysis would then continue as before but with $a_1 p_1+a_2 p_2$ as a joint wavenumber.
In this situation however the wavenumbers have to be handled differently. For example, we can have separately
large values of $|p_1|$ and $|p_2|$ but the value of $a_1p_1+a_2p_2$ can be small.

\section{Conclusions}

The numerical analysis of hypoelliptic PDEs with variable coefficients and Dirac initial datum is notoriously difficult because standard approaches to convergence analysis are not directly  applicable. In contrast with parabolic initial-value problems, in the case of hypoelliptic PDEs the situation
is complicated by the fact that
diffusion generally acts only in some, but not all, co-ordinate directions. The approach of \cite{GC} for the one-dimensional heat equation is based on Fourier analysis to show approximation of a certain order for the low wavenumber components and exponential decay of the high wavenumber components for smoothing schemes. In the present case of hypoelliptic equations the analysis is more intricate because of the interplay of wavenumbers in the different co-ordinate directions as a consequence of variable coefficients and the resulting noncommutativity of the spatial differential operators in the different directions.

We exploit the special structure of the Kolmogorov equations under consideration by performing a Fourier transform in the direction with the pure transport term, and discretize in the diffusive direction(s) by a finite difference scheme. This results in a parametrized system of semi-discrete equations, which can be solved by standard methods, and the resulting solution is then transformed back from Fourier space using exponentially convergent numerical quadrature, thus leading to an effective technique for dimensional reduction. The numerical analysis sheds light on the interaction of different wavenumbers and allows us to derive second-order convergence for the proposed numerical scheme in the absence of time-discretization.

An open question at present is how to replicate the analysis herein for the fully discrete counterpart of the method, which also includes time-discretization. Although unconditional stability of the backward Euler scheme for initial data in $L^2$, can be proved by standard energy estimates, a similar argument does not seem to apply in the case of Dirac initial datum, as the use of a discrete negative-order Sobolev norm, which would be the natural choice of norm for the discrete representation of the Dirac function on the finite difference grid, in conjunction with discrete counterparts of parabolic smoothing estimates in the diffusive directions, is obstructed by
the particular form of the variable coefficients in the Kolmogorov equations under consideration, and the complex nature of error propagation for the fully-discrete scheme has not so far allowed a rigorous convergence analysis of the fully discrete scheme.

\bibliographystyle{plain}
\bibliography{fourier}



\appendix

\section{Implementation of the FFT}
\label{app:ifft}

Starting from our inversion formula, for a given $i$ and $m$, we have
\[  U_{j,k,m}\sim \frac{\Delta p}{2\pi}    \sum^{l_{max}}_{l=l_{min}}   V_{j,l,m}\, {\rm e}^{-ip_ly_k}\]
for a suitable range of $l$-values, and we choose $l_{min}=-l_{max}$. We write $Y_k=U_{j,k,m}$ and $X_l=V_{j,l,m}$, so that
\[  Y_k\sim \frac{\Delta p}{2\pi}    \sum^{l_{max}}_{l=l_{min}}   X_l \,{\rm e}^{-ip_ly_k}.\]
We then have
\begin{align*}
\frac{\Delta p}{2\pi}    \sum^{l_{max}}_{l=l_{min}}   X_l\,{\rm e}^{-ip_ly_k}&=\frac{\Delta p}{2\pi}    \sum^{l_{max}-l_{min}+1}_{\tilde{l}=1}
 X_{\tilde{l} + l_{min}-1}\,{\rm e}^{-i\Delta p(\tilde{l} + l_{min}-1)\Delta yk}\\
&=\frac{\Delta p}{2\pi}    \sum^{l_{max}-l_{min}+1}_{\tilde{l}=1}   X_{\tilde{l} + l_{min}-1}\,{\rm e}^{-\frac{2\pi i}{N}\frac{N\Delta p\Delta y}{2\pi}(\tilde{l} + l_{min}-1)k}.
\end{align*}
To align this with the MATLAB procedure, we require
%
%
%
\[  N=\frac{2\pi}{\Delta p\Delta y}      \]
and $N=l_{max}-l_{min}+1$. We then have
\begin{align*}
\frac{\Delta p}{2\pi}    \sum^{N}_{l=1}   X_{l+ l_{min}-1}\,{\rm e}^{-\frac{2\pi i}{N}(l+ l_{min}-1)k}&=\frac{\Delta p}{2\pi}    \sum^{N}_{l=1}    X_{l+ l_{min}-1}
\,{\rm e}^{-\frac{2\pi i}{N}(l-1)k}\,{\rm e}^{-\frac{2\pi i}{N} l_{min}k}\\
&=\frac{1}{N}\sum^{N}_{l=1} \left({\rm e}^{-\frac{2\pi i}{N} l_{min}( \tilde{k}   -1    )}\frac{N\Delta p}{2\pi}    \right)   X_{l+ l_{min}-1}\,{\rm e}^{-\frac{2\pi i}{N}(l-1)( \tilde{k}   -1    )}
\end{align*}
so that
\begin{align}
Y_{ \tilde{k}   -1  }\sim\frac{1}{N}\sum^{N}_{l=1} W_{l , \tilde{k} }\,{\rm e}^{-\frac{2\pi i}{N}(l-1)( \tilde{k}   -1    )}=IFFT(W),
\end{align}
where, for each value of $\tilde{k}$,
\[ W_{l , \tilde{k} } =  \left({\rm e}^{-\frac{2\pi i}{N} l_{min}( \tilde{k}   -1    )}\frac{N\Delta p}{2\pi}    \right) X_{l+ l_{min}-1} .   \]
This expresses our inversion in terms of the MATLAB IFFT. 
%
%
If the $p$ grid values are $(  l_{min}\Delta p,\ldots, l_{max}\Delta p         )$, then we have
\[p_{range}=p_{max}-p_{min}=\Delta p(l_{max}-l_{min})\]
and
\[ \Delta y=\frac{2\pi}{N\Delta p }=\frac{2\pi}{(l_{max}-l_{min}+1)\Delta p}  =\frac{2\pi}{ p_{range}+\Delta p  } \sim  \frac{2\pi}{p_{range}}  .  \]
Due to this definition, care has to be taken when applying these procedures to evaluate the Fourier transform
above. If we fix $p_{range}$ and then refine the $p$-grid, the $y$-grid resolution will not change. If we want to provide values of the approximation at more
closely-spaced $y$-values, we will have to increase the $p$-range.

\end{document}

%% file: macpap2.tex
\newcommand{\be}{\begin{equation}}
\newcommand{\ee}{\end{equation}}
\newcommand{\bean}{\begin{eqnarray}}
\newcommand{\eean}{\end{eqnarray}}
\newcommand{\bea}{\begin{eqnarray*}}
\newcommand{\eea}{\end{eqnarray*}}
\newcommand{\bfig}{\begin{figure}}
\newcommand{\efig}{\end{figure}}
\newcommand{\ba}{\begin{array}}
\newcommand{\ea}{\end{array}}
